\documentclass[preprint,12pt]{elsarticle}

\usepackage{amssymb}
\usepackage{subfigure}
\usepackage{graphicx}
\usepackage{epstopdf}
 \usepackage{amsthm}
\usepackage{latexsym}
\usepackage{amsmath}
\usepackage{amssymb}
\usepackage{diagbox}
\usepackage{slashbox}
\usepackage{color}
\usepackage{arydshln}

\newcommand{\udots}{\mathinner{\mskip1mu\raise1pt\vbox{\kern7pt\hbox{.}}
\mskip2mu\raise4pt\hbox{.}\mskip2mu\raise7pt\hbox{.}\mskip1mu}}
\usepackage{amsmath}
\newtheorem{definition}{Definition}
\newtheorem{theorem}{Theorem}

\begin{document}

\begin{frontmatter}



\title{ Multivariate double truncated expectation and covariance risk measures for elliptical distributions
}

\author{Baishuai  Zuo}
\author{Chuancun Yin\corref{cor1}}
\cortext[cor1]{Corresponding author.}
 \ead{ccyin@qfnu.edu.cn}

\address{School of Statistics, Qufu Normal University, Qufu, Shandong 273165, P. R. China}

\begin{abstract}
 The main objective of this work is to calculate the multivariate double truncated expectation (MDTE) and covariance (MDTCov) for elliptical distributions. We also consider double truncated expectation (DTE) and variance (DTV) for univariate elliptical distributions. The exact expressions of MDTE and MDTCov are derived for some special cases of the family, such as normal, student-$t$, logistic, Laplace and Pearson type VII distributions. As numerical illustration, the DTE, DTV, MDTE and MDTCov for normal distribution are computed in details. Finally, we discuss MDTE and MDTCov of three industry segments' (Banks, Insurance, Financial and Credit Service)  stock return in London stock exchange.
\end{abstract}


\end{frontmatter}

\baselineskip=20pt

\section{Introduction}
 In the last decade, much research has been devoted to risk measures and many multi-dimensional extensions have been investigated. For example, several generalizations of the classical univariate conditional tail expectation have been proposed (see, e.g., Cousin and Bernardino, 2014; Cai et al., 2017; Shushi and Yao, 2020; Ahmed et al., 2021; Cai et al., 2021; Ortega-Jim$\acute{e}$nez et al., 2021). Specially,
Landsman et al. (2016) defined a new multivariate tail conditional expectation (MTCE):
\begin{align}\label{(3)}
\nonumber\mathrm{MTCE}_{\boldsymbol{q}}(\mathbf{X})&=\mathrm{E}\left[\mathbf{X}|\mathbf{X}>VaR_{\boldsymbol{q}}(\mathbf{X})\right]\\
&=\mathrm{E}[\mathbf{X}|X_{1}>VaR_{q_{1}}(X_{1}),\cdots,X_{n}>VaR_{q_{n}}(X_{n})],
\end{align}
 where
 $\boldsymbol{q}=(q_{1},\cdots,q_{n})\in(0,~1)^{n},$
 $\mathbf{X}=(X_{1},~X_{2},\cdots,X_{n})^{T}$ is an $n\times1$ vector of risks with cumulative distribution function (cdf) $F_{\mathbf{X}}(\boldsymbol{x})$ and tail function $\overline{F}_{\mathbf{X}}(\boldsymbol{x})$, $$VaR_{\boldsymbol{q}}(\mathbf{X})=(VaR_{q_{1}}(X_{1}),~VaR_{q_{2}}(X_{2}),\cdots,VaR_{q_{n}}(X_{n}))^{T},$$ and $VaR_{q_{k}}(X_{k}),~k=1,~2,\cdots,n$ is the value at risk (VaR) measure of the random variable $X_{k}$, being the $q_{k}$-th quantile of $X_{k}$. The MTCE is reduced to tail conditional expectation (TCE) when $n=1$. Since Lansman et al. (2016) has been derived formula of MTCE for elliptical distribution, Mousavi et al. (2019) derived expression of MTCE for scale
mixtures of skew-normal distribution, Zuo and Yin (2021a) and Zuo and Yin (2021b) extended those results to generalized skew-elliptical and location-scale
mixture of elliptical distributions, respectively.

Landsman et al. (2018) defined a novel type of a multivariate tail covariance (MTCov):
\begin{align}\label{(4)}
\nonumber\mathrm{MTCov}_{\boldsymbol{q}}(\mathbf{X})&=\mathrm{E}\left[(\mathbf{X}-\mathrm{MTCE}_{\boldsymbol{q}}(\mathbf{X}))(\mathbf{X}-\mathrm{MTCE}_{\boldsymbol{q}}(\mathbf{X}))^{T}|\mathbf{X}>VaR_{\boldsymbol{q}}(\mathbf{X})\right]\\
&=\inf_{\boldsymbol{c}\in\mathbb{R}^{n}}\mathrm{E}\left[(\mathbf{X}-\boldsymbol{c})(\mathbf{X}-\boldsymbol{c})^{T}|\mathbf{X}>VaR_{\boldsymbol{q}}(\mathbf{X})\right],
\end{align}
which is an extension of tail variance (TV) measure.
Especially, Landsman et al. (2018) studied MTCov for elliptical distribution, Landsman et al. (2021) extended this result to logistic elliptical distributions, and Zuo and Yin (2021c)  extended it to generalized skew-elliptical distributions.

 In Arismendi and Broda (2017), the authors derived multivariate elliptical truncated moment generating function, and first second-order moments of quadratic forms of the multivariate normal, Student and generalized hyperbolic distributions. They pointed
out that elliptical truncated moments' expansions was applied to many areas, such as the design of experiment (see Thompson, 1976), robust estimation (see Cuesta-Albertos et al., 2008), outlier detections (see, e.g., Riani et al., 2009; Cerioli, 2010), robust regression (see Torti et al., 2012), robust detection (see Cerioli et al. (2014)), statistical estates' estimation (see Shi et al., 2014) and risk averse selection (see Hanasusanto et al., 2015). Furthermore, in Ogasawara (2021), the author derived a non-recursive formula for various moments of the
multivariate normal distribution with sectional truncation, and introduced the importance of truncated moments in biological field, such as animals or plants breeding programs (see Herrend\"{o}rfer and Tuchscherer, 1996) and medical treatments with risk variables as blood pressures and pulses, where low and high values of the variables are of primary concern.

 In Roozegar et al. (2020), the authors derived explicit expressions of the first two moments for doubly truncated multivariate normal
mean-variance mixture distributions. Inspired by those work, we introduced  two multivariate risk measures: multivariate double truncated expectation and covariance, and derived MDTE and MDTCov for elliptical distributions. We also presented expressions of MDTE and MDTCov for some special cases of this distributions, including normal, student-$t$, logistic, Laplace and Pearson type VII distributions.

The rest of the paper is organized as follows. Section 2 introduces the definitions and properties of multivariate double truncated expectation and covariance risk measures. Section 3 introduces elliptical family, including normal, student-$t$, logistic, Laplace and Pearson type VII distributions. In Section 4, we derive MDTE for elliptical distributions, and some special cases are presented in this section. In Section 5, we also derive MDTCov for elliptical distributions, and some special cases are shown in this section. We give numerical illustration and illustrative example in Sections 6 and 7, respectively. Finally, Section 8,  is the concluding remarks.
\section{MDTE and MDTCov}
In this section, we introduce two multivariate risk measures. They are multivariate double truncated expectation (MDTE) and  multivariate double truncated covariance (MDTCov), which are generalizations of multivariate tail conditional expectation (MTCE) and multivariate tail covariance matrix (MTCov), respectively.
\begin{definition}\label{def.1}
For an $n\times1$ vector $\mathbf{X}$, multivariate double truncated expectation $(\mathrm{MDTE})$ is defined by
\begin{align}\label{(5)}
\nonumber\mathrm{MDTE}_{(\boldsymbol{p},\boldsymbol{q})}(\mathbf{X})&=\mathrm{E}\left[\mathbf{X}|VaR_{\boldsymbol{p}}(\mathbf{X})<\mathbf{X}<VaR_{\boldsymbol{q}}(\mathbf{X})\right]\\
&=\mathrm{E}[\mathbf{X}|VaR_{p_{i}}(X_{i})<X_{i}<VaR_{q_{i}}(X_{i}), i=1,2,\cdots,n],
\end{align}
where  $$\boldsymbol{p}=(p_{1},\cdots,p_{n}),~\boldsymbol{q}=(q_{1},\cdots,q_{n})\in(0,~1)^{n},~p_{k}<q_{k},~k=1,2,\cdots,n,$$
$$VaR_{\boldsymbol{v}}(\mathbf{X})=(VaR_{v_{1}}(X_{1}),~VaR_{v_{2}}(X_{2}),\cdots,VaR_{v_{n}}(X_{n}))^{T},~\boldsymbol{v}=\boldsymbol{p,q},$$ and $VaR_{v_{k}}(X_{k}),~k=1,~2,\cdots,n$ is the value at risk (VaR) measure of the random variable $X_{k}$, being the $v_{k}$-th quantile of $X_{k}$.
\end{definition}
\noindent $\mathbf{Remark~1}$ When $n=1$, the multivariate double truncated expectation $(\mathrm{MDTE})$ risk measure is reduced to double truncated expectation $(\mathrm{DTE})$:
\begin{align}\label{(6)}
&\mathrm{DTE}_{(p,q)}(X)=\mathrm{E}\left[X|VaR_{p}(X)<X<VaR_{q}(X)\right],~p,~q\in(0,~1), ~and~p<q.
\end{align}
\begin{definition}\label{def.2}
For an $n\times1$ vector $\mathbf{X}$, multivariate double truncated covariance $\mathrm{(MDTCov)}$ is defined by
\begin{align}\label{(7)}
\nonumber&\mathrm{MDTCov}_{(\boldsymbol{p},\boldsymbol{q})}(\mathbf{X})=\\
&\mathrm{E}\left[(\mathbf{X}-\mathrm{MDTE}_{(\boldsymbol{p},\boldsymbol{q})}(\mathbf{X}))(\mathbf{X}-\mathrm{MDTE}_{(\boldsymbol{p},\boldsymbol{q})}(\mathbf{X}))^{T}|VaR_{\boldsymbol{p}}(\mathbf{X})<\mathbf{X}<VaR_{\boldsymbol{q}}(\mathbf{X})\right].
\end{align}
\end{definition}
\noindent$\mathbf{Remark~2}$ When $n=1$, the multivariate double truncated covariance\\
 $\mathrm{(MDTCov)}$ risk measure is reduced to double truncated variance $\mathrm{(DTV)}$:
\begin{align}\label{(8)}
\nonumber&\mathrm{DTV}_{(p,q)}(X)\\
&=\mathrm{E}\left[(X-\mathrm{DTE}_{(p,q)}(X))^{2}|VaR_{p}(X)<X<VaR_{q}(X)\right].
\end{align}

Note that the MDTE risk measure is MTCE as $VaR_{\boldsymbol{q}}(\mathbf{X})\rightarrow\boldsymbol{+\infty}$, and the MDTCov risk measure is MTCov as $VaR_{\boldsymbol{q}}(\mathbf{X})\rightarrow\boldsymbol{+\infty}$. The MDTE and MDTCov risk measures are multivariate expectation and multivariate covariance matix as $VaR_{\boldsymbol{p}}(\mathbf{X})\rightarrow\boldsymbol{-\infty}$ and  $VaR_{\boldsymbol{q}}(\mathbf{X})\rightarrow\boldsymbol{+\infty}$, respectively. When $n=1$, DTE and DTV risk measure are TCE and TV as $VaR_{q}(X)\rightarrow +\infty$, respectively. In addition, the MDTCov is different from following multivariate double truncated conditional covariance (MDTCCov):
\begin{align}\label{(9)}
\nonumber&\mathrm{MDTCCov}_{(\boldsymbol{p},\boldsymbol{q})}(\mathbf{X})\\
\nonumber&=\mathrm{E}\left[(\mathbf{X}-\mathrm{E}(\mathbf{X}))(\mathbf{X}-\mathrm{E}(\mathbf{X}))^{T}|VaR_{\boldsymbol{p}}(\mathbf{X})<\mathbf{X}<VaR_{\boldsymbol{q}}(\mathbf{X})\right]\\
\nonumber&=\mathrm{MDTCov}_{(\boldsymbol{p},\boldsymbol{q})}(\mathbf{X})+\mathrm{MDTE}_{(\boldsymbol{p},\boldsymbol{q})}(\mathbf{X})\mathrm{MDTE}^{T}_{(\boldsymbol{p},\boldsymbol{q})}(\mathbf{X})\\
&~~~-\mathrm{MDTE}_{(\boldsymbol{p},\boldsymbol{q})}(\mathbf{X})\mathrm{E}^{T}(\mathbf{X})-\mathrm{E}(\mathbf{X})\mathrm{MDTE}^{T}_{(\boldsymbol{p},\boldsymbol{q})}(\mathbf{X})+\mathrm{E}(\mathbf{X})\mathrm{E}^{T}(\mathbf{X}).
\end{align}

The following are some properties of the MDTE and MDTCov risk measures. The proof is trivial.\\
\noindent$\mathbf{Proposition~1}$ For any $n\times1$ random vectors $\mathbf{X}=(X_{1},X_{2},\cdots,X_{n})^{T}$ and $\mathbf{Y}=(Y_{1},Y_{2},\cdots,Y_{n})^{T}$, the MDTE risk measure has following properties:\\
(i) For any positive constant $b$, we have
$$\mathrm{MDTE}_{(\boldsymbol{p},\boldsymbol{q})}(b\mathbf{X})=b\mathrm{MDTE}_{(\boldsymbol{p},\boldsymbol{q})}(\mathbf{X});$$
(ii) For any vector of constants $\boldsymbol{\gamma}\in \mathbb{R}^{n}$
$$\mathrm{MDTE}_{(\boldsymbol{p},\boldsymbol{q})}(\mathbf{X}+\boldsymbol{\gamma})=\mathrm{MDTE}_{(\boldsymbol{p},\boldsymbol{q})}(\mathbf{X})+\boldsymbol{\gamma};$$
(iii) When $\mathbf{X}$ has independent components, we have
$$\mathrm{MDTE}_{(\boldsymbol{p},\boldsymbol{q})}(\mathbf{X})=(\mathrm{DTE}_{(p_{1},q_{1})}(X_{1}),\mathrm{DTE}_{(p_{2},q_{2})}(X_{2}),\cdots,\mathrm{DTE}_{(p_{n},q_{n})}(X_{n}))^T;$$
(iv) If $\mathbf{Y}\stackrel{a.s.}\geq\mathbf{X}$, then
$$\mathrm{MDTE}_{(\boldsymbol{p},\boldsymbol{q})}(\mathbf{Y}-\mathbf{X})\geq \boldsymbol{0},$$
where $\boldsymbol{0}$ is vector of $n$ zeros.\\
(v) If $\mathbf{S}=(\mathbf{X}^T,\mathbf{Y}^T)^T$ is elliptically distributed vector, then
$$\mathrm{MDTE}_{(\boldsymbol{p},\boldsymbol{q})}(\mathbf{Y}+\mathbf{X})\leq\mathrm{MDTE}_{(\boldsymbol{p},\boldsymbol{q})}(\mathbf{S})_{1}+\mathrm{MDTE}_{(\boldsymbol{p},\boldsymbol{q})}(\mathbf{S})_{2},$$
where $\mathrm{MDTE}_{(\boldsymbol{p},\boldsymbol{q})}(\mathbf{S})_{1}$ and $\mathrm{MDTE}_{(\boldsymbol{p},\boldsymbol{q})}(\mathbf{S})_{2}$ are vectors of the first $n$ elements and the last $n$ elements of $\mathrm{MDTE}_{(\boldsymbol{p},\boldsymbol{q})}(\mathbf{S})$, respectively. When $\mathbf{X}$ and $\mathbf{Y}$ are independent, then the equality holds.\\
\noindent$\mathbf{Proposition~2}$ For any $n\times1$ random vector $\mathbf{X}=(X_{1},X_{2},\cdots,X_{n})^{T}$, the MDTCov risk measure has following properties:\\
(i) For any positive constant $b$, we have
$$\mathrm{MDTCov}_{(\boldsymbol{p},\boldsymbol{q})}(b\mathbf{X})=b^{2}\mathrm{MDTCov}_{(\boldsymbol{p},\boldsymbol{q})}(\mathbf{X});$$
(ii) For any vector of constants $\boldsymbol{\gamma}\in \mathbb{R}^{n}$
$$\mathrm{MDTCov}_{(\boldsymbol{p},\boldsymbol{q})}(\mathbf{X}+\boldsymbol{\gamma})=\mathrm{MDTCov}_{(\boldsymbol{p},\boldsymbol{q})}(\mathbf{X});$$
(iii) When $\mathbf{X}$ has independent components, we have
$$\mathrm{MDTCov}_{(\boldsymbol{p},\boldsymbol{q})}(\mathbf{X})=diag(\mathrm{DTV}_{(p_{1},q_{1})}(X_{1}),\mathrm{DTV}_{(p_{2},q_{2})}(X_{2}),\cdots,\mathrm{DTV}_{(p_{n},q_{n})}(X_{n})),$$
where $diag()$ is diagonal matrix.

Similarly, from matrix MDTCov we can obtain the multivariate double truncated correlation (MDTCorr)
matrix:
\begin{align}\label{(a9)}
\mathrm{MDTCorr}_{(\boldsymbol{p},\boldsymbol{q})}(\mathbf{X})=\bigg(\frac{\mathrm{MDTCov}_{(\boldsymbol{p},\boldsymbol{q})}(\mathbf{X})_{ij}}{\sqrt{\mathrm{MDTCov}_{(\boldsymbol{p},\boldsymbol{q})}(\mathbf{X})_{ii}}\sqrt{\mathrm{MDTCov}_{(\boldsymbol{p},\boldsymbol{q})}(\mathbf{X})_{jj}}}\bigg)_{ij=1,\cdots,n}.
\end{align}

\section{Family of elliptical distributions}
 An $n \times 1$ random vector $X = (X_1,\cdots, X_n)^{T}$ is said to have an elliptically
symmetric distribution if it's probability density function (pdf) exists, the form will be (see Landsman and Valdez, 2003)
\begin{align}\label{(10)}
f_{\boldsymbol{X}}(\boldsymbol{x}):=\frac{c_{n}}{\sqrt{|\boldsymbol{\Sigma}|}}g_{n}\left\{\frac{1}{2}(\boldsymbol{x}-\boldsymbol{\mu})^{T}\mathbf{\Sigma}^{-1}(\boldsymbol{x}-\boldsymbol{\mu})\right\},~\boldsymbol{x}\in\mathbb{R}^{n},
\end{align}
 where $\boldsymbol{\mu}$ is an $n\times1$ location vector, $\mathbf{\Sigma}$ is an $n\times n$ scale matrix, and $g_{n}(u)$, $u\geq0$, is the density generator of $\mathbf{X}$. We denote it by $\mathbf{X}\sim E_{n}(\boldsymbol{\mu},\boldsymbol{\Sigma},g_{n})$. The density generator $g_{n}$ satisfies the condition
\begin{align}\label{(11)}
\int_{0}^{\infty}s^{n/2-1}g_{n}(s)\mathrm{d}s<\infty,
\end{align}
 and the normalizing constant  $c_n$ is given by
\begin{align*}
c_{n}=\frac{\Gamma(n/2)}{(2\pi)^{n/2}}\left[\int_{0}^{\infty}s^{n/2-1}g_{n}(s)\mathrm{d}s\right]^{-1}.
\end{align*}

Cumulative generator $\overline{G}_{n}(u)$ and $\overline{\mathcal{G}}_{n}(u)$ are defined as follows:
\begin{align*}
\overline{G}_{n}(u)=\int_{u}^{\infty}{g}_{n}(v)\mathrm{d}v
\end{align*}
and
\begin{align*}
\overline{\mathcal{G}}_{n}(u)=\int_{u}^{\infty}{G}_{n}(v)\mathrm{d}v,
\end{align*}
and their normalizing constants are, respectively, written:
\begin{align*}
c_{n}^{\ast}=\frac{\Gamma(n/2)}{(2\pi)^{n/2}}\left[\int_{0}^{\infty}s^{n/2-1}\overline{G}_{n}(s)\mathrm{d}s\right]^{-1}
\end{align*}
and
\begin{align*}
c_{n}^{\ast\ast}=\frac{\Gamma(n/2)}{(2\pi)^{n/2}}\left[\int_{0}^{\infty}s^{n/2-1}\overline{\mathcal{G}}_{n}(s)\mathrm{d}s\right]^{-1}.
\end{align*}
$\mathbf{X}^{\ast}\sim E_{n}(\boldsymbol{\mu},~\boldsymbol{\Sigma},~\overline{G}_{n})$ and $\mathbf{X}^{\ast\ast}\sim E_{n}(\boldsymbol{\mu},~\boldsymbol{\Sigma},~\overline{\mathcal{G}}_{n})$ (see Zuo et al., 2021) are respectively called elliptical random vectors with generators $\overline{G}_{n}(u)$ and $\overline{\mathcal{G}}_{n}(u)$, if their density functions (if them exist) defined by
\begin{align}\label{(12)}
f_{\boldsymbol{X}^{\ast}}(\boldsymbol{x})=\frac{c_{n}^{\ast}}{\sqrt{|\boldsymbol{\Sigma}|}}\overline{G}_{n}\left\{\frac{1}{2}(\boldsymbol{x}-\boldsymbol{\mu})^{T}\mathbf{\Sigma}^{-1}(\boldsymbol{x}-\boldsymbol{\mu})\right\},~\boldsymbol{x}\in\mathbb{R}^{n}
\end{align}
\begin{align}\label{(13)}
f_{\boldsymbol{X}^{\ast\ast}}(\boldsymbol{x})=\frac{c_{n}^{\ast\ast}}{\sqrt{|\boldsymbol{\Sigma}|}}\overline{\mathcal{G}}_{n}\left\{\frac{1}{2}(\boldsymbol{x}-\boldsymbol{\mu})^{T}\mathbf{\Sigma}^{-1}(\boldsymbol{x}-\boldsymbol{\mu})\right\},~\boldsymbol{x}\in\mathbb{R}^{n}.
\end{align}
Here the density generators $\overline{G}_{n}(u)$ and $\overline{\mathcal{G}}_{n}(u)$ satisfy the conditions:
\begin{align}\label{(14)}
\int_{0}^{\infty}s^{n/2-1}\overline{G}_{n}(s)\mathrm{d}s<+\infty
\end{align}
and
\begin{align}\label{(15)}
\int_{0}^{\infty}s^{n/2-1}\overline{\mathcal{G}}_{n}(s)\mathrm{d}s<+\infty.
\end{align}

The following are some special members of the class of elliptical distributions.\\
$\mathbf{Example~1}$ (Multivariate normal distribution). Suppose that $\mathbf{X}\sim N_{n}\left(\boldsymbol{\mu},~\boldsymbol{\Sigma}\right)$. In this case,
 the density generators   are  
\begin{align}\label{(a17)}
g(u)=\overline{G}(u)=\overline{\mathcal{G}}(u)=\exp\{-u\},
\end{align}
and the normalizing constants are written:
\begin{align}\label{(b17)}
c_{n}=c_{n}^{\ast}=c_{n}^{\ast\ast}=(2\pi)^{-\frac{n}{2}}.
\end{align}
 $\mathbf{Example~2}$ (Multivariate student-$t$ distribution). Suppose that
 \begin{align}\label{(18)}
 \mathbf{X}\sim St_{n}\left(\boldsymbol{\mu},~\boldsymbol{\Sigma},~m\right).
 \end{align}
 In this case, the density generators $g_{n}(u)$, $\overline{G}_{n}(u)$ and $\overline{\mathcal{G}}_{n}(u)$ are expressed (for details see Zuo et al., 2021):
\begin{align}\label{(a18)}
g_{n}(u)=\left(1+\frac{2u}{m}\right)^{-(m+n)/2},
\end{align}
\begin{align}\label{(b18)}
\overline{G}_{n}(u)=\frac{m}{m+n-2}\left(1+\frac{2u}{m}\right)^{-(m+n-2)/2}
\end{align}
and
\begin{align}\label{(c18)}
\overline{\mathcal{G}}_{n}(u)=\frac{m}{m+n-2}\frac{m}{m+n-4}\left(1+\frac{2u}{m}\right)^{-(m+n-4)/2}.
\end{align}
The normalizing constants are written:
\begin{align}\label{(d18)}
 c_{n}=\frac{\Gamma\left((m+n)/2\right)}{\Gamma(m/2)(m\pi)^{\frac{n}{2}}},
 \end{align}
\begin{align}\label{(e18)}
 \nonumber c_{n}^{\ast}&=\frac{(m+n-2)\Gamma(n/2)}{(2\pi)^{n/2}m}\left[\int_{0}^{\infty}u^{n/2-1}\left(1+\frac{2t}{m}\right)^{-(m+n-2)/2}\mathrm{d}u\right]^{-1}\\
 &=\frac{(m+n-2)\Gamma(n/2)}{(m\pi)^{n/2}mB(\frac{n}{2},~\frac{m-2}{2})},~if~m>2
 \end{align} and
 \begin{align}\label{(f18)}
 \nonumber c_{n}^{\ast\ast}&=\frac{(m+n-2)(m+n-4)\Gamma(n/2)}{(2\pi)^{n/2}m^{2}}\left[\int_{0}^{\infty}u^{n/2-1}\left(1+\frac{2t}{m}\right)^{-(m+n-4)/2}\mathrm{d}u\right]^{-1}\\
 &=\frac{(m+n-2)(m+n-4)\Gamma(n/2)}{(m\pi)^{n/2}m^{2}B(\frac{n}{2},~\frac{m-4}{2})},~if~m>4,
 \end{align}
 where $\Gamma(\cdot)$ and $B(\cdot,\cdot)$ are Gamma function and Beta function, respectively. \\
$\mathbf{Example~3}$ (Multivariate logistic distribution). Suppose that $\mathbf{X}\sim Lo_{n}\left(\boldsymbol{\mu},~\boldsymbol{\Sigma}\right)$. In this case, the density generators $g_{n}(u)$, $\overline{G}_{n}(u)$ and $\overline{\mathcal{G}}_{n}(u)$ are expressed (for details see Zuo et al., 2021):
 \begin{align}\label{(a19)}
g_{n}(u)=\frac{\exp(-u)}{[1+\exp(-u)]^{2}},
\end{align}
\begin{align}\label{(b19)}
\overline{G}_{n}(u)=\frac{\exp(-u)}{1+\exp(-u)}
\end{align}
and
\begin{align}\label{(c19)}
 \overline{\mathcal{G}}_{n}(u)=\ln\left[1+\exp(-u)\right].
 \end{align}
 The normalizing constants are written:
 \begin{align}\label{(d19)}
\nonumber c_{n}&=\frac{\Gamma(n/2)}{(2\pi)^{n/2}}\left[\int_{0}^{\infty}u^{n/2-1}\frac{\exp(-u)}{[1+\exp(-u)]^{2}}\mathrm{d}u\right]^{-1}\\
&=\frac{1}{(2\pi)^{n/2}\Psi_{2}^{\ast}(-1,\frac{n}{2},1)},
\end{align}
\begin{align}\label{(e19)}
 \nonumber c_{n}^{\ast}&=\frac{\Gamma(n/2)}{(2\pi)^{n/2}}\left[\int_{0}^{\infty}u^{n/2-1}\frac{\exp(-u)}{1+\exp(-u)}\mathrm{d}u\right]^{-1}\\
 &=\frac{1}{(2\pi)^{n/2}\Psi_{1}^{\ast}(-1,\frac{n}{2},1)}
\end{align}
 and
 \begin{align}\label{(f19)}
 \nonumber c_{n}^{\ast\ast}&=\frac{\Gamma(n/2)}{(2\pi)^{n/2}}\left\{\int_{0}^{\infty}u^{n/2-1}\ln\left[1+\exp(-u)\right]\mathrm{d}u\right\}^{-1}\\
 \nonumber&=\frac{\Gamma(n/2)}{(2\pi)^{n/2}}\left[\frac{2}{n}\int_{0}^{\infty}u^{n/2}\frac{e^{-u}}{1+e^{-u}}\mathrm{d}u\right]^{-1}\\
 &=\frac{1}{(2\pi)^{n/2}\Psi_{1}^{\ast}(-1,\frac{n}{2}+1,1)}.
\end{align}
$\mathbf{Remark~3}$ Here $\Psi_{\kappa}^{\ast}(z,s,a)$ is the generalized Hurwitz-Lerch zeta function defined by (see
Lin et al., 2006)
$$\Psi_{\kappa}^{\ast}(z,s,a)=\frac{1}{\Gamma(\kappa)}\sum_{n=0}^{\infty}\frac{\Gamma(\kappa+n)}{n!}\frac{z^{n}}{(n+a)^{s}},$$
which has an integral representation
$$\Psi_{\kappa}^{\ast}(z,s,a)=\frac{1}{\Gamma(s)}\int_{0}^{\infty}\frac{t^{s-1}e^{-at}}{(1-ze^{-t})^{\kappa}}\mathrm{d}t,$$
where $\mathcal{R}(a)>0$, $\mathcal{R}(s)>0$ when $|z|\leq1~(z\neq1)$, $\mathcal{R}(s)>1$ when $z=1$.\\
 $\mathbf{Example~4}$ (Multivariate Laplace distribution). Suppose  $\mathbf{X}\sim La_{n}\left(\boldsymbol{\mu},~\boldsymbol{\Sigma}\right)$. In this case, the density generators $g_{n}(u)$, $\overline{G}_{n}(u)$ and $\overline{\mathcal{G}}_{n}(u)$ are expressed (for details see Zuo et al., 2021):
\begin{align}\label{(a20)}
g_{n}(u)=\exp(-\sqrt{2u}),
\end{align}
\begin{align}\label{(b20)}
\overline{G}_{n}(u)=(1+\sqrt{2u})\exp(-\sqrt{2u})
\end{align}
and
\begin{align}\label{(c20)}
 \overline{\mathcal{G}}_{n}(u)=(3+2u+3\sqrt{2u})\exp(-\sqrt{2u}).
 \end{align}
The normalizing constants are written:
\begin{align}\label{(d20)}
c_{n}=\frac{\Gamma(n/2)}{2\pi^{n/2}\Gamma(n)},
 \end{align}
\begin{align}\label{(e20)}
c_{n}^{\ast}=\frac{n\Gamma(n/2)}{2\pi^{n/2}\Gamma(n+2)}
\end{align}
and
\begin{align}\label{(f20)}
 c_{n}^{\ast\ast}=\frac{n(n+2)\Gamma(n/2)}{2\pi^{n/2}\Gamma(n+4)}.
 \end{align}
$\mathbf{Example~5}$ (Multivariate Pearson type VII distribution). Suppose that
$$
\mathbf{X}\sim PVII_{n}\left(\boldsymbol{\mu},~\boldsymbol{\Sigma},~t\right).
$$
In this case, the density generators $g_{n}(u)$, $\overline{G}_{n}(u)$ and $\overline{\mathcal{G}}_{n}(u)$ are expressed:
 \begin{align}\label{(a21)}
 g_{n}(u)=(1+2u)^{-t},
 \end{align}
 \begin{align}\label{(b21)}
 \overline{G}_{n}(u)=\frac{1}{2(t-1)}(1+2u)^{-(t-1)}
 \end{align}
 and
 \begin{align}\label{(c21)}
 \overline{\mathcal{G}}_{n}(u)=\frac{1}{4(t-1)(t-2)}(1+2u)^{-(t-2)}.
 \end{align}
The normalizing constants are written:
\begin{align}\label{(d21)}
 c_{n}=\frac{\Gamma\left(t\right)}{\Gamma(t-n/2)\pi^{\frac{n}{2}}},~t>\frac{n}{2},
 \end{align}
\begin{align}\label{(e21)}
 c_{n}^{\ast}=\frac{\Gamma(\frac{n}{2})2(t-1)}{\pi^{\frac{n}{2}}B(\frac{n}{2},t-1-\frac{n}{2})},~t>1+\frac{n}{2}
 \end{align}
 and
 \begin{align}\label{(f21)}
 c_{n}^{\ast\ast}=\frac{\Gamma(\frac{n}{2})4(t-1)(t-2)}{\pi^{\frac{n}{2}}B(\frac{n}{2},t-2-\frac{n}{2})},~t>2+\frac{n}{2}.
 \end{align}
 $\mathbf{Remark~4}$ Multivariate Pearson type VII distribution is related to the multivariate student-t by the transformation $\mathbf{X}=[\mathbf{Y}-(1-\sqrt{m}\boldsymbol{\mu})]/\sqrt{m}$, where $\mathbf{Y}$ is as in (\ref{(18)}) with $m=2t-n$ (see Zografos, 2008).
\section{Multivariate double truncated expectation}
Consider a random vector  $\mathbf{X}\sim E_{n}\left(\boldsymbol{\mu},~\boldsymbol{\Sigma},~g_{n}\right)$ with finite vector $\boldsymbol{\mu}=(\mu_1,\cdots,\mu_n)^{T}$, positive defined matrix ${\bf {\Sigma}}= (\sigma_{ij})_{i,j=1}^{n}$ and probability density function $f_{\boldsymbol{X}}(\boldsymbol{x})$.

Let $\mathbf{Y}=\mathbf{\Sigma}^{-\frac{1}{2}}(\mathbf{X}-\boldsymbol{\mu})\sim E_{n}\left(\boldsymbol{0},~\boldsymbol{I_{n}},~g_{n}\right).$
Writing
$$\boldsymbol{\eta_{v}}=\left(\eta_{\boldsymbol{v},1},~\eta_{\boldsymbol{v},2},\cdots,\eta_{\boldsymbol{v},n}\right)^{T}=\mathbf{\Sigma}^{-\frac{1}{2}}(\boldsymbol{x_{v}-\mu}),$$
where
$$\boldsymbol{x_{v}}=VaR_{\boldsymbol{v}}(\boldsymbol{X}),~ \boldsymbol{\eta}_{\boldsymbol{v},-k}=\left(\eta_{\boldsymbol{v},1},~\eta_{\boldsymbol{v},2},\cdots,\eta_{\boldsymbol{v},k-1},~\eta_{\boldsymbol{v},k+1},\cdots,\eta_{\boldsymbol{v},n}\right)^{T}$$ and $$\boldsymbol{\eta}_{\boldsymbol{v},-k,j}=\left(\eta_{\boldsymbol{v},1},\cdots,\eta_{\boldsymbol{v},k-1},~\eta_{\boldsymbol{v},k+1},\cdots,\eta_{\boldsymbol{v},j-1},~\eta_{\boldsymbol{v},j+1},\cdots,\eta_{\boldsymbol{v},n}\right)^{T}, ~\boldsymbol{v}=\boldsymbol{p},\boldsymbol{q}.$$
To present MDTE,
we define an new truncated distribution function as follows: $$F_{\mathbf{Z}}(\boldsymbol{a},\boldsymbol{b})=\int_{\boldsymbol{a}}^{\boldsymbol{b}}f_{\mathbf{Z}}(\boldsymbol{z})\mathrm{d}\boldsymbol{z},$$
where $f_{\mathbf{Z}}(\boldsymbol{z})$ is pdf of random vector $\mathbf{Z}$.
We find that $F_{\mathbf{Z}}(\boldsymbol{a},\boldsymbol{b})=F_{\mathbf{Z}}(\boldsymbol{b})$ as $\boldsymbol{a}\rightarrow\mathbf{-\infty}$, and $F_{\mathbf{Z}}(\boldsymbol{a},\boldsymbol{b})=\overline{F}_{\mathbf{Z}}(\boldsymbol{b})$ as $\boldsymbol{b}\rightarrow\mathbf{+\infty}$.\\
We define shifted cumulative generator
\begin{align}\label{(16)}
\overline{G}_{n-1,\boldsymbol{v},k}^{\ast}(u)=\overline{G}_{n}(u+\frac{1}{2}\eta_{\boldsymbol{v},k}^{2}),~k=1,2,\cdots,n,~\boldsymbol{v}=\boldsymbol{p},\boldsymbol{q},
\end{align}
and normalizing constant
\begin{align}\label{(a16)}
c_{n-1,\boldsymbol{v},k}^{\ast}=\frac{\Gamma(\frac{n-1}{2})}{(2\pi)^{(n-1)/2}}\left[\int_{0}^{\infty}s^{(n-3)/2}\overline{G}^{\ast}_{n-1,\boldsymbol{v},k}(s)\mathrm{d}s\right]^{-1}.
\end{align}
\begin{theorem}\label{th.2}  Let $\mathbf{X}\sim E_{n}(\boldsymbol{\mu},~\mathbf{\Sigma},~g_{n})$ $(n\geq2)$ be an $n$-dimensional elliptical random vector with density generator $g_{n}$, positive defined scale matrix ${\bf {\Sigma}}$ and finite expectation $\boldsymbol{\mu}$. Further, it satisfies conditions (\ref{(11)}) and (\ref{(14)}).
 Then,
\begin{align}\label{(22)} &\mathrm{MDTE}_{(\boldsymbol{p,q})}(\mathbf{X})=\boldsymbol{\mu}+\mathbf{\Sigma}^{\frac{1}{2}}\frac{\boldsymbol{\delta_{\boldsymbol{p,q}}}}{F_{\mathbf{Y}}(\boldsymbol{\eta}_{\boldsymbol{p}},\boldsymbol{\eta}_{\boldsymbol{q}})},
\end{align}
where
$$\boldsymbol{\delta_{\boldsymbol{p,q}}}=\left(\delta_{1,\boldsymbol{p,q}},~\delta_{2,\boldsymbol{p,q}},\cdots,\delta_{n,\boldsymbol{p,q}}\right)^{T}$$
with
\begin{align*}
&\delta_{k,\boldsymbol{p,q}}\\
=&\frac{c_{n}}{c_{n-1,\boldsymbol{p},k}^{\ast}}F_{\mathbf{Y}_{\boldsymbol{p},-k}}(\boldsymbol{\eta}_{\boldsymbol{p},-k},\boldsymbol{\eta}_{\boldsymbol{q},-k})-\frac{c_{n}}{c_{n-1,\boldsymbol{q},k}^{\ast}}F_{\mathbf{Y}_{\boldsymbol{q},-k}}(\boldsymbol{\eta}_{\boldsymbol{p},-k},\boldsymbol{\eta}_{\boldsymbol{q},-k}),~k=1,~2,\cdots,n,
 \end{align*}
 $\mathbf{Y}\sim E_{n}\left(\boldsymbol{0},~\boldsymbol{I_{n}},~g_{n}\right)$, $\mathbf{Y}_{\boldsymbol{v},-k}\sim E_{n-1}\left(\boldsymbol{0},~\boldsymbol{I_{n-1}},~\overline{G}_{n-1,\boldsymbol{v},k}^{\ast}\right),$
 and $\overline{G}_{n-1,\boldsymbol{v},k}^{\ast}(u)$ and
$c_{n-1,\boldsymbol{v},k}^{\ast}$ are same as defined in (\ref{(16)}) and (\ref{(a16)}), respectively.
\end{theorem}
\noindent Proof. Using definition, we have
\begin{align*} &\mathrm{MDTE}_{(\boldsymbol{p,q})}(\mathbf{X})=\frac{c_{n}}{\sqrt{|\mathbf{\Sigma}|}F_{\mathbf{X}}(\boldsymbol{x_{p},x_{q}})}\int_{\boldsymbol{x_{p}}}^{\boldsymbol{x_{q}}}\boldsymbol{x}g_{n}\left(\frac{1}{2}\boldsymbol{(x-\mu)^{T}\Sigma^{-1}(x-\mu)}\right)\mathrm{d}\boldsymbol{x}.
\end{align*}
Applying translation $\mathbf{Y}=\mathbf{\Sigma}^{-\frac{1}{2}}(\mathbf{X}-\boldsymbol{\mu})$, we obtain
\begin{align*} \mathrm{MDTE}_{(\boldsymbol{p,q})}(\mathbf{X})&=\frac{c_{n}}{F_{\mathbf{Y}}(\boldsymbol{\eta_{p},\eta_{q}})}\int_{\boldsymbol{\eta_{p}}}^{\boldsymbol{\eta_{q}}}\left(\mathbf{\Sigma}^{\frac{1}{2}}\boldsymbol{y}+\boldsymbol{\mu}\right) g_{n}\left(\frac{1}{2}\boldsymbol{y^{T}y}\right)\mathrm{d}\boldsymbol{y}\\
&=\boldsymbol{\mu}+\frac{\mathbf{\Sigma}^{\frac{1}{2}}}{F_{\mathbf{Y}}(\boldsymbol{\eta_{p},\eta_{q}})}\boldsymbol{\delta_{p,q}},
\end{align*}
where $$\boldsymbol{\delta_{p,q}}=\int_{\boldsymbol{\eta_{p}}}^{\boldsymbol{\eta_{q}}}c_{n}\boldsymbol{y} g_{n}\left(\frac{1}{2}\boldsymbol{y^{T}y}\right)\mathrm{d}\boldsymbol{y}.$$
Note that
\begin{align*}
\delta_{k,\boldsymbol{p,q}}&=c_{n}\int_{\boldsymbol{\eta_{p}}}^{\boldsymbol{\eta_{q}}}y_{k}g_{n}\left(\frac{1}{2}\boldsymbol{y^{T}y}\right)\mathrm{d}\boldsymbol{y}\\
&=c_{n}\int_{\boldsymbol{\eta}_{\boldsymbol{p},-k}}^{\boldsymbol{\eta}_{\boldsymbol{q},-k}}\int_{\eta_{\boldsymbol{p},k}}^{\eta_{\boldsymbol{q},k}}-\partial_{k}\overline{G}_{n}\left(\frac{1}{2}\boldsymbol{y}_{-k}^{T}\boldsymbol{y}_{-k}+\frac{1}{2}y_{k}^{2}\right)\mathrm{d}\boldsymbol{y}_{-k}\\
&=c_{n}\int_{\boldsymbol{\eta}_{\boldsymbol{p},-k}}^{\boldsymbol{\eta}_{\boldsymbol{q},-k}}\left[\overline{G}_{n}\left(\frac{1}{2}\boldsymbol{y}_{-k}^{T}\boldsymbol{y}_{-k}+\frac{1}{2}\eta_{\boldsymbol{p},k}^{2}\right)-\overline{G}_{n}\left(\frac{1}{2}\boldsymbol{y}_{-k}^{T}\boldsymbol{y}_{-k}+\frac{1}{2}\eta_{\boldsymbol{q},k}^{2}\right)\right]\mathrm{d}\boldsymbol{y}_{-k}\\
&=c_{n}\left\{\frac{1}{c_{n-1,\boldsymbol{p},k}^{\ast}}F_{\mathbf{Y}_{\boldsymbol{p},-k}}(\boldsymbol{\eta}_{\boldsymbol{p},-k},\boldsymbol{\eta}_{\boldsymbol{q},-k})-\frac{1}{c_{n-1,\boldsymbol{q},k}^{\ast}}F_{\mathbf{Y}_{\boldsymbol{q},-k}}(\boldsymbol{\eta}_{\boldsymbol{p},-k},\boldsymbol{\eta}_{\boldsymbol{q},-k})\right\},
\end{align*}
$k=1,~2,\cdots,n.$ Therefore, we obtain the desired result.

When $n=1$, the DTE for univariate elliptical distribution is given:
\begin{align}\label{(23)}
&\mathrm{DTE}_{(p,q)}(X)=\mu+\sigma\frac{\delta_{p,q}}{F_{Y}(\eta_{p},\eta_{q})},
\end{align}
where
\begin{align*}
\delta_{p,q}=c_{1}\left[\overline{G}_{1}\left(\frac{1}{2}\eta_{p}^{2}\right)-\overline{G}_{1}\left(\frac{1}{2}\eta_{q}^{2}\right)\right],
\end{align*}
and $Y\sim E_{1}\left(0,~1,~g_{1}\right)$.\\
$\mathbf{Remark~5}$ Note that $\boldsymbol{x_{q}}=VaR_{\boldsymbol{q}}(\mathbf{X})\rightarrow\boldsymbol{+\infty}$ in Theorem $1$, one gets the formula of Theorem 1 in Landsman et al. (2018), which is generalization of Theorem 1 in Landsman et al. (2016); When $x_{q}=VaR_{q}(X)\rightarrow +\infty$ in (\ref{(23)}), we obtain result of Theorem 1 in  Landsman and Valdez (2003), which is the generalization of Theorem 1 in Jiang et al. (2016).\\
\noindent$\mathbf{Corollary~1}$
  Let $\mathbf{X}\sim N_{n}(\boldsymbol{\mu},~\mathbf{\Sigma})$ $(n\geq2)$. Then
\begin{align}\label{(24)} &\mathrm{MDTE}_{(\boldsymbol{p,q})}(\mathbf{X})=\boldsymbol{\mu}+\mathbf{\Sigma}^{\frac{1}{2}}\frac{\boldsymbol{\delta_{\boldsymbol{p,q}}}}{F_{\mathbf{Y}}(\boldsymbol{\eta}_{\boldsymbol{p}},\boldsymbol{\eta}_{\boldsymbol{q}})},
\end{align}
where
$$\boldsymbol{\delta_{\boldsymbol{p,q}}}=\left(\delta_{1,\boldsymbol{p,q}},~\delta_{2,\boldsymbol{p,q}},\cdots,\delta_{n,\boldsymbol{p,q}}\right)^{T}$$
with
\begin{align*}
\delta_{k,\boldsymbol{p,q}}&=\left[\phi(\eta_{\boldsymbol{p},k})-\phi(\eta_{\boldsymbol{q},k})\right]F_{\mathbf{Y}_{-k}}(\boldsymbol{\eta}_{\boldsymbol{p},-k},\boldsymbol{\eta}_{\boldsymbol{q},-k}),~k=1,~2,\cdots,n,
 \end{align*}
$\mathbf{Y}\sim N_{n}\left(\boldsymbol{0},~\boldsymbol{I_{n}}\right)$, $\mathbf{Y}_{-k}\sim N_{n-1}\left(\boldsymbol{0},~\boldsymbol{I_{n-1}}\right)$, and $\phi(\cdot)$ is pdf of $1$-dimensional standard normal distribution. \\
When $n=1$, we obtain the DTE for univariate normal distribution:
\begin{align}\label{(25)}
&\mathrm{DTE}_{(p,q)}(X)=\mu+\sigma\frac{\delta_{p,q}}{F_{Y}(\eta_{p},\eta_{q})},
\end{align}
where
\begin{align}\label{(a25)}
\delta_{p,q}=\phi(\eta_{p})-\phi(\eta_{q}),
\end{align}
and $Y\sim N_{1}\left(0,~1\right)$.\\
\noindent$\mathbf{Corollary~2}$
  Let $\mathbf{X}\sim St_{n}(\boldsymbol{\mu},~\mathbf{\Sigma},~m)$ $(n\geq2)$. Then
\begin{align}\label{(26)} &\mathrm{MDTE}_{(\boldsymbol{p,q})}(\mathbf{X})=\boldsymbol{\mu}+\mathbf{\Sigma}^{\frac{1}{2}}\frac{\boldsymbol{\delta_{\boldsymbol{p,q}}}}{F_{\mathbf{Y}}(\boldsymbol{\eta}_{\boldsymbol{p}},\boldsymbol{\eta}_{\boldsymbol{q}})},
\end{align}
where
$$\boldsymbol{\delta_{\boldsymbol{p,q}}}=\left(\delta_{1,\boldsymbol{p,q}},~\delta_{2,\boldsymbol{p,q}},\cdots,\delta_{n,\boldsymbol{p,q}}\right)^{T}$$
with
\begin{align*}
&\delta_{k,\boldsymbol{p,q}}=\\
&\frac{c_{n}}{c_{n-1,\boldsymbol{p},k}^{\ast}}F_{\mathbf{Y}_{\boldsymbol{p},-k}}(\boldsymbol{\eta}_{\boldsymbol{p},-k},\boldsymbol{\eta}_{\boldsymbol{q},-k})-\frac{c_{n}}{c_{n-1,\boldsymbol{q},k}^{\ast}}F_{\mathbf{Y}_{\boldsymbol{q},-k}}(\boldsymbol{\eta}_{\boldsymbol{p},-k},\boldsymbol{\eta}_{\boldsymbol{q},-k}),~k=1,~2,\cdots,n,
 \end{align*}
 $\mathbf{Y}\sim St_{n}\left(\boldsymbol{0},~\boldsymbol{I_{n}},~m\right)$, $\mathbf{Y}_{\boldsymbol{v},-k}\sim St_{n-1}\left(\boldsymbol{0},~\Delta_{\boldsymbol{v},k},~m-1\right),$
 $$\Delta_{\boldsymbol{v},k}=\left[\frac{m\left(1+\eta_{\boldsymbol{v},k}^{2}/m\right)}{m-1}\right]\mathbf{I}_{n-1} $$ and
$$c_{n-1,\boldsymbol{v},k}^{\ast}=\frac{\Gamma\left(\frac{m+n-2}{2}\right)(m+n-2)}{\Gamma\left(\frac{m-1}{2}\right)[(m-1)\pi]^{(n-1)/2}m}\left(1+\frac{\eta_{\boldsymbol{v},k}^{2}}{m}\right)^{(m+n-2)/2},~\boldsymbol{v}=\boldsymbol{p},\boldsymbol{q}.$$
We can further simplify
$$\frac{c_{n}}{c_{\boldsymbol{v},k}^{\ast}}=\frac{\Gamma\left(\frac{m-1}{2}\right)\left(\frac{m-1}{m}\right)^{(n-1)/2}}{\Gamma\left(\frac{m}{2}\right)\sqrt{\pi/m}}\left(1+\frac{\eta_{\boldsymbol{v},k}^{2}}{m}\right)^{-(m+n-2)/2},~\boldsymbol{v}=\boldsymbol{p},\boldsymbol{q}.$$
When $n=1$, the DTE of univariate student-t distribution is given:
\begin{align}\label{(27)}
&\mathrm{DTE}_{(p,q)}(X)=\mu+\sigma\frac{\delta_{p,q}}{F_{Y}(\eta_{p},\eta_{q})},
\end{align}
where
\begin{align}\label{(a27)}
\delta_{p,q}=\frac{\Gamma((m+1)/2)m}{\Gamma(m/2)\sqrt{m\pi}(m-1)}\left[\left(1+\frac{\eta_{p}^{2}}{m}\right)^{-(m-1)/2}-\left(1+\frac{\eta_{q}^{2}}{m}\right)^{-(m-1)/2}\right],
\end{align}
and $Y\sim St_{1}\left(0,~1,~m\right)$.\\
\noindent$\mathbf{Corollary~3}$ Let $\mathbf{X}\sim Lo_{n}(\boldsymbol{\mu},~\mathbf{\Sigma})$ $(n\geq2)$. Then
\begin{align}\label{(28)} &\mathrm{MDTE}_{(\boldsymbol{p,q})}(\mathbf{X})=\boldsymbol{\mu}+\mathbf{\Sigma}^{\frac{1}{2}}\frac{\boldsymbol{\delta_{\boldsymbol{p,q}}}}{F_{\mathbf{Y}}(\boldsymbol{\eta}_{\boldsymbol{p}},\boldsymbol{\eta}_{\boldsymbol{q}})},
\end{align}
where
$$\boldsymbol{\delta_{\boldsymbol{p,q}}}=\left(\delta_{1,\boldsymbol{p,q}},~\delta_{2,\boldsymbol{p,q}},\cdots,\delta_{n,\boldsymbol{p,q}}\right)^{T}$$
with
\begin{align*}
&\delta_{k,\boldsymbol{p,q}}=\\
&\frac{c_{n}}{c_{n-1,\boldsymbol{p},k}^{\ast}}F_{\mathbf{Y}_{\boldsymbol{p},-k}}(\boldsymbol{\eta}_{\boldsymbol{p},-k},\boldsymbol{\eta}_{\boldsymbol{q},-k})-\frac{c_{n}}{c_{n-1,\boldsymbol{q},k}^{\ast}}F_{\mathbf{Y}_{\boldsymbol{q},-k}}(\boldsymbol{\eta}_{\boldsymbol{p},-k},\boldsymbol{\eta}_{\boldsymbol{q},-k}),~k=1,~2,\cdots,n,
 \end{align*}
 $\mathbf{Y}\sim Lo_{n}\left(\boldsymbol{0},~\boldsymbol{I_{n}}\right)$,
\begin{align*}
c_{n-1,\boldsymbol{v},k}^{\ast}&=\frac{\Gamma((n-1)/2)\exp\left\{\frac{\eta_{\boldsymbol{v},k}^{2}}{2}\right\}}{(2\pi)^{(n-1)/2}}\left[\int_{0}^{\infty}\frac{t^{(n-3)/2}\exp\{-t\}}{1+\exp\left\{-\frac{\eta_{\boldsymbol{v},k}^{2}}{2}\right\}\exp\{-t\}}\mathrm{d}t\right]^{-1}\\
&=\frac{\exp\left\{\frac{\eta_{\boldsymbol{v},k}^{2}}{2}\right\}}{(2\pi)^{(n-1)/2}\Psi_{1}^{\ast}\left(-\sqrt{2\pi}\phi(\eta_{\boldsymbol{v},k}),~\frac{n-1}{2},~1\right)},
\end{align*}and pdf of $\mathbf{Y}_{\boldsymbol{v},-k}$:
 \begin{align*}
&f_{\mathbf{Y}_{\boldsymbol{v},-k}}(\boldsymbol{t})\\
&=c_{n-1,\boldsymbol{v},k}^{\ast}\frac{\exp\left\{-\frac{\boldsymbol{t}^{T}\boldsymbol{t}}{2}-\frac{\eta_{\boldsymbol{v},k}^{2}}{2}\right\}}{1+\exp\left\{-\frac{\boldsymbol{t}^{T}\boldsymbol{t}}{2}-\frac{\eta_{\boldsymbol{v},k}^{2}}{2}\right\}}, ~\boldsymbol{t}\in\mathbb{R}^{n-1},~\boldsymbol{v=p},\boldsymbol{q},~k=1,~2,\cdots,n.
\end{align*}
We can further simplify
$$\frac{c_{n}}{c_{n-1,\boldsymbol{v},k}^{\ast}}=\frac{\Psi_{1}^{\ast}\left(-\sqrt{2\pi}\phi(\eta_{\boldsymbol{v},k}),~\frac{n-1}{2},~1\right)\phi(\eta_{\boldsymbol{v},k})}{\Psi_{2}^{\ast}\left(-1,\frac{n}{2},1\right)},~k=1,2,\cdots,n,~\boldsymbol{v=p},\boldsymbol{q}.$$
When $n=1$,we obtain the DTE of univariate logistic distribution:
\begin{align}\label{(29)}
&\mathrm{DTE}_{(p,q)}(X)=\mu+\sigma\frac{\delta_{p,q}}{F_{Y}(\eta_{p},\eta_{q})},
\end{align}
where
\begin{align}\label{(a29)}
\delta_{p,q}=\frac{1}{\Psi^{\ast}_{2}\left(-1,\frac{1}{2},1\right)}\left[\frac{\phi(\eta_{p})}{1+\sqrt{2\pi}\phi(\eta_{p})}-\frac{\phi(\eta_{q})}{1+\sqrt{2\pi}\phi(\eta_{q})}\right],
\end{align}
and $Y\sim Lo_{1}\left(0,~1\right)$.\\
\noindent$\mathbf{Corollary~4}$
  Let $\mathbf{X}\sim La_{n}(\boldsymbol{\mu},~\mathbf{\Sigma})$ $(n\geq2)$. Then
\begin{align}\label{(30)} &\mathrm{MDTE}_{(\boldsymbol{p,q})}(\mathbf{X})=\boldsymbol{\mu}+\mathbf{\Sigma}^{\frac{1}{2}}\frac{\boldsymbol{\delta_{\boldsymbol{p,q}}}}{F_{\mathbf{Y}}(\boldsymbol{\eta}_{\boldsymbol{p}},\boldsymbol{\eta}_{\boldsymbol{q}})},
\end{align}
where
$$\boldsymbol{\delta_{\boldsymbol{p,q}}}=\left(\delta_{1,\boldsymbol{p,q}},~\delta_{2,\boldsymbol{p,q}},\cdots,\delta_{n,\boldsymbol{p,q}}\right)^{T}$$
with
\begin{align*}
&\delta_{k,\boldsymbol{p,q}}=\\
&\frac{c_{n}}{c_{n-1,\boldsymbol{p},k}^{\ast}}F_{\mathbf{Y}_{\boldsymbol{p},-k}}(\boldsymbol{\eta}_{\boldsymbol{p},-k},\boldsymbol{\eta}_{\boldsymbol{q},-k})-\frac{c_{n}}{c_{n-1,\boldsymbol{q},k}^{\ast}}F_{\mathbf{Y}_{\boldsymbol{q},-k}}(\boldsymbol{\eta}_{\boldsymbol{p},-k},\boldsymbol{\eta}_{\boldsymbol{q},-k}),~k=1,~2,\cdots,n,
 \end{align*}
 $\mathbf{Y}\sim La_{n}\left(\boldsymbol{0},~\boldsymbol{I_{n}}\right)$,
$$c_{n-1,\boldsymbol{v},k}^{\ast}=\frac{\Gamma\left(\frac{n-1}{2}\right)}{(2\pi)^{(n-1)/2}}\left[\int_{0}^{\infty}t^{\frac{n-3}{2}}\left(1+\sqrt{2t+\eta_{\boldsymbol{v}}^{2},k}\right)\exp\left\{-\sqrt{2t+\eta_{\boldsymbol{v},k}^{2}}\right\}\mathrm{d}t\right]^{-1},$$
and pdf of $\mathbf{Y}_{\boldsymbol{v},-k}$:
\begin{align*}
f_{\mathbf{Y}_{\boldsymbol{v},-k}}(\boldsymbol{t})=c_{n-1,\boldsymbol{v},k}^{\ast}\left(1+\sqrt{\boldsymbol{t}^{T}\boldsymbol{t}+\eta_{\boldsymbol{v},k}^{2}}\right)\exp\left\{-\sqrt{\boldsymbol{t}^{T}\boldsymbol{t}+\eta_{\boldsymbol{v},k}^{2}}\right\},~k=1,~2,\cdots,n,
\end{align*}
$\boldsymbol{t}\in\mathbb{R}^{n-1}$,~$\boldsymbol{v}=\boldsymbol{p},\boldsymbol{q}$.\\
We can further simplify
$$\frac{c_{n}}{c_{n-1,\boldsymbol{v},k}^{\ast}}=\frac{\Gamma\left(\frac{n}{2}\right)2^{(n-3)/2}}{\sqrt{\pi}\Gamma\left(\frac{n-1}{2}\right)\Gamma(n)}\left[\int_{0}^{\infty}t^{\frac{n-3}{2}}\left(1+\sqrt{2t+\eta_{\boldsymbol{v}}^{2},k}\right)\exp\left\{-\sqrt{2t+\eta_{\boldsymbol{v},k}^{2}}\right\}\mathrm{d}t\right],$$
$~k=1,~2,\cdots,n,~\boldsymbol{v}=\boldsymbol{p},\boldsymbol{q}.$\\
When $n=1$, the DTE of univariate Laplace distribution is given:
\begin{align}\label{(31)}
&\mathrm{DTE}_{(p,q)}(X)=\mu+\sigma\frac{\delta_{p,q}}{F_{Y}(\eta_{p},\eta_{q})},
\end{align}
where
\begin{align}\label{(a31)}
\delta_{p,q}=\frac{1}{2}\left[(1+|\eta_{p}|)\exp(-|\eta_{p}|)-(1+|\eta_{q}|)\exp(-|\eta_{q}|)\right],
\end{align}
$Y\sim La_{1}\left(0,~1\right)$, and $|\cdot|$ is the absolute value function.\\
\noindent$\mathbf{Corollary~5}$
  Let $\mathbf{X}\sim PVII_{n}(\boldsymbol{\mu},~\mathbf{\Sigma},~t)$ $(n\geq2)$. Then
\begin{align}\label{(32)} &\mathrm{MDTE}_{(\boldsymbol{p,q})}(\mathbf{X})=\boldsymbol{\mu}+\mathbf{\Sigma}^{\frac{1}{2}}\frac{\boldsymbol{\delta_{\boldsymbol{p,q}}}}{F_{\mathbf{Y}}(\boldsymbol{\eta}_{\boldsymbol{p}},\boldsymbol{\eta}_{\boldsymbol{q}})},
\end{align}
where
$$\boldsymbol{\delta_{\boldsymbol{p,q}}}=\left(\delta_{1,\boldsymbol{p,q}},~\delta_{2,\boldsymbol{p,q}},\cdots,\delta_{n,\boldsymbol{p,q}}\right)^{T}$$
with
\begin{align*}
&\delta_{k,\boldsymbol{p,q}}=\\
&\frac{c_{n}}{c_{n-1,\boldsymbol{p},k}^{\ast}}F_{\mathbf{Y}_{\boldsymbol{p},-k}}(\boldsymbol{\eta}_{\boldsymbol{p},-k},\boldsymbol{\eta}_{\boldsymbol{q},-k})-\frac{c_{n}}{c_{n-1,\boldsymbol{q},k}^{\ast}}F_{\mathbf{Y}_{\boldsymbol{q},-k}}(\boldsymbol{\eta}_{\boldsymbol{p},-k},\boldsymbol{\eta}_{\boldsymbol{q},-k}),~k=1,~2,\cdots,n,
 \end{align*}
 $\mathbf{Y}\sim PVII_{n}\left(\boldsymbol{0},~\boldsymbol{I_{n}},~t\right)$, $\mathbf{Y}_{\boldsymbol{v},-k}\sim PVII_{n-1}\left(\boldsymbol{0},~\Lambda_{\boldsymbol{v},k},~t-1\right),$
 $$\Lambda_{\boldsymbol{v},k}=\left(1+\eta_{\boldsymbol{v},k}^{2}\right)\mathbf{I}_{n-1} $$ and
$$c_{n-1,\boldsymbol{v},k}^{\ast}=\frac{2\Gamma(t)}{\Gamma\left(t-\frac{n-1}{2}\right)\pi^{(n-1)/2}}\left(1+\eta_{\boldsymbol{v},k}^{2}\right)^{t-1},~\boldsymbol{v}=\boldsymbol{p},\boldsymbol{q}.$$
We can further simplify
$$\frac{c_{n}}{c_{n-1,\boldsymbol{v},k}^{\ast}}=\frac{\Gamma\left(t-\frac{n-1}{2}\right)}{2\Gamma\left(t-\frac{n}{2}\right)\sqrt{\pi}}\left(1+\eta_{\boldsymbol{v},k}^{2}\right)^{-(t-1)},~\boldsymbol{v}=\boldsymbol{p},\boldsymbol{q}.$$
When $n=1$, we obtain the DTE of univariate Pearson type VII distribution:
\begin{align}\label{(33)}
&\mathrm{DTE}_{(p,q)}(X)=\mu+\sigma\frac{\delta_{p,q}}{F_{Y}(\eta_{p},\eta_{q})},
\end{align}
where
\begin{align}\label{(a33)}
\delta_{p,q}=\frac{\Gamma(t-1)}{2\Gamma(t-1/2)\sqrt{\pi}}\left[\left(1+\eta_{p}^{2}\right)^{-(t-1)}-\left(1+\eta_{q}^{2}\right)^{-(t-1)}\right],
\end{align}
and $Y\sim PVII_{1}\left(0,~1,~t\right)$.
 \section{Multivariate double truncated covariance}
To derive formula for MDTCov we also define shifted cumulative generator
\begin{align}\label{(34)}
\overline{\mathcal{G}}_{n-2,\boldsymbol{u}r,\boldsymbol{v}l,-ij}^{\ast\ast}(t)=\overline{\mathcal{G}}_{n}\left(t+\frac{1}{2}\eta_{\boldsymbol{u},r}^{2}+\frac{1}{2}\eta_{\boldsymbol{v},l}^{2}\right),~\boldsymbol{u}=\boldsymbol{p},\boldsymbol{q};\boldsymbol{v}=\boldsymbol{p},\boldsymbol{q};r=i,j;l=i,j, \end{align}
and normalizing constant
\begin{align}\label{(35)}
c_{n-2,\boldsymbol{u}r,\boldsymbol{v}l,ij}^{\ast\ast}=\frac{\Gamma(\frac{n-2}{2})}{(2\pi)^{(n-2)/2}}\left[\int_{0}^{\infty}s^{(n-4)/2}\overline{\mathcal{G}}_{n-2,\boldsymbol{u}r,\boldsymbol{v}l,-ij}^{\ast\ast}(t)\mathrm{d}t\right]^{-1}.
\end{align}

In the following, we formulate the theorem that gives multivariate double truncated covariance (MDTCov) for elliptical distributions.
\begin{theorem}\label{th.2} Under conditions (\ref{(11)}), (\ref{(14)}) and (\ref{(15)}),
 $\mathrm{MDTCov}$ of elliptical distributed random vector $\mathbf{X}\sim E_{n}\left(\boldsymbol{\mu},~\mathbf{\Sigma},~g_{n}\right)$ is give by
\begin{align}\label{(36)} &\mathrm{MDTCov}_{(\boldsymbol{p},\boldsymbol{q})}(\mathbf{X})=\mathbf{\Sigma}^{\frac{1}{2}}\mathbf{\Upsilon}_{\boldsymbol{p},\boldsymbol{q}}\mathbf{\Sigma}^{\frac{1}{2}},
\end{align}
where
\begin{align*}
&\mathbf{\Upsilon}_{\boldsymbol{p,q},i,j}=\frac{1}{F_{\mathbf{Y}}(\boldsymbol{\eta_{p}},\boldsymbol{\eta_{q}})}\\
& \boldsymbol{\cdot}\bigg\{\frac{c_{n}}{c_{n-2,\boldsymbol{p}i,\boldsymbol{p}j,ij}^{\ast\ast}}F_{\mathbf{Y}_{\boldsymbol{p}i,\boldsymbol{p}j,-ij}}(\boldsymbol{\eta}_{\boldsymbol{p},-ij},\boldsymbol{\eta}_{\boldsymbol{q},-ij})-\frac{c_{n}}{c_{n-2,\boldsymbol{p}i,\boldsymbol{q}j,ij}^{\ast\ast}}F_{\mathbf{Y}_{\boldsymbol{p}i,\boldsymbol{q}j,-ij}}(\boldsymbol{\eta}_{\boldsymbol{p},-ij},\boldsymbol{\eta}_{\boldsymbol{q},-ij})\\
&+\frac{c_{n}}{c_{n-2,\boldsymbol{q}i,\boldsymbol{q}j,ij}^{\ast\ast}}F_{\mathbf{Y}_{\boldsymbol{q}i,\boldsymbol{q}j,-ij}}(\boldsymbol{\eta}_{\boldsymbol{p},-ij},\boldsymbol{\eta}_{\boldsymbol{q},-ij}) -\frac{c_{n}}{c_{n-2,\boldsymbol{p}j,\boldsymbol{q}i,ij}^{\ast\ast}}F_{\mathbf{Y}_{\boldsymbol{p}j,\boldsymbol{q}i,-ij}}(\boldsymbol{\eta}_{\boldsymbol{p},-ij},\boldsymbol{\eta}_{\boldsymbol{q},-ij}) \bigg\}\\
&-\mathrm{MDTE}_{(\boldsymbol{\eta_{p}},\boldsymbol{\eta_{q}})}(\mathbf{Y})_{i}\mathrm{MDTE}_{(\boldsymbol{\eta_{p}},\boldsymbol{\eta_{q}})}(\mathbf{Y})_{j}
,~i\neq j,~i=1,2,\cdots,n,~j=1,2,\cdots,n
\end{align*}
and
\begin{align*}
&\mathbf{\Upsilon}_{\boldsymbol{p,q},ii}=\\ &\frac{1}{F_{\mathbf{Y}}(\boldsymbol{\eta}_{\boldsymbol{p}},\boldsymbol{\eta}_{\boldsymbol{q}})}\bigg\{\frac{c_{n}}{c_{n-1,\boldsymbol{p},i}^{\ast}}\eta_{\boldsymbol{p},i}F_{\mathbf{Y}_{\boldsymbol{p},-i}}(\boldsymbol{\eta}_{\boldsymbol{p},-i},\boldsymbol{\eta}_{\boldsymbol{q},-i}) -\frac{c_{n}}{c_{n-1,\boldsymbol{q},i}^{\ast}}\eta_{\boldsymbol{q},i}F_{\mathbf{Y}_{\boldsymbol{q},-i}}(\boldsymbol{\eta}_{\boldsymbol{p},-i},\boldsymbol{\eta}_{\boldsymbol{q},-i}) \\ &~~~+\frac{c_{n}}{c_{n}^{\ast}}F_{\mathbf{Y}^{\ast}}(\boldsymbol{\eta}_{\boldsymbol{p}},\boldsymbol{\eta}_{\boldsymbol{q}})\bigg\}-[\mathrm{MDTE}_{(\boldsymbol{\eta_{p}},\boldsymbol{\eta_{q}})}(\mathbf{Y})_{i}]^{2},~i=1,2,\cdots,n,
\end{align*}
with
\begin{align}\label{(37)}
\nonumber&\mathrm{MDTE}_{(\boldsymbol{\eta_{p}},\boldsymbol{\eta_{q}})}(\mathbf{Y})_{k}\\
&=\frac{1}{F_{\mathbf{Y}}(\boldsymbol{\eta}_{\boldsymbol{p}},\boldsymbol{\eta}_{\boldsymbol{q}})}\bigg\{\frac{c_{n}}{c_{n-1,\boldsymbol{p},k}^{\ast}}F_{\mathbf{Y}_{\boldsymbol{p},-k}}(\boldsymbol{\eta}_{\boldsymbol{p},-k},\boldsymbol{\eta}_{\boldsymbol{q},-k})-\frac{c_{n}}{c_{n-1,\boldsymbol{q},k}^{\ast}}F_{\mathbf{Y}_{\boldsymbol{q},-k}}(\boldsymbol{\eta}_{\boldsymbol{p},-k},\boldsymbol{\eta}_{\boldsymbol{q},-k})\bigg\},
\end{align}
  $k=1,~2,\cdots,n,$ $\mathbf{Y}\sim E_{n}\left(\boldsymbol{0},~\boldsymbol{I_{n}},~g_{n}\right)$, $\mathbf{Y}^{\ast}\sim E_{n}\left(\boldsymbol{0},~\boldsymbol{I_{n}},~\overline{G}_{n}\right)$,\\ $\mathbf{Y}_{\boldsymbol{u}r,\boldsymbol{v}l,-ij}\sim E_{n-2}(\boldsymbol{0},~\mathbf{I_{n-2}},~\overline{\mathcal{G}}_{n-2,\boldsymbol{u}r,\boldsymbol{v}l,-ij}^{\ast\ast})$, $\boldsymbol{u}=\boldsymbol{p},\boldsymbol{q};\boldsymbol{v}=\boldsymbol{p},\boldsymbol{q};r=i,j;l=i,j,$ $\mathbf{Y}_{\boldsymbol{v},-k}\sim E_{n-1}\left(\boldsymbol{0},~\boldsymbol{I_{n-1}},~\overline{G}_{n-1,\boldsymbol{v},k}^{\ast}\right),$
 and $\overline{G}_{n-1,\boldsymbol{v},k}^{\ast}(u)$, $\overline{\mathcal{G}}_{n-2,\boldsymbol{u}r,\boldsymbol{v}l,-ij}^{\ast\ast}$,
$c_{n-1,\boldsymbol{v},k}^{\ast}$ and $c_{n-2,\boldsymbol{u}r,\boldsymbol{v}l,ij}^{\ast\ast}$ are same as defined in (\ref{(16)}), (\ref{(34)}), (\ref{(a16)}) and (\ref{(35)}), respectively.
\end{theorem}
\noindent Proof. Using definition of MDTCov, we have
\begin{align*}
&\mathrm{MDTCov}_{\boldsymbol{p},\boldsymbol{q}}(\mathbf{X})\\
&=\mathrm{E}\left[(\mathbf{X}-\mathrm{MDTE}_{\boldsymbol{p},\boldsymbol{q}}(\mathbf{X}))(\mathbf{X}-\mathrm{MDTE}_{\boldsymbol{p},\boldsymbol{q}}(\mathbf{X}))^{T}|VaR_{\boldsymbol{p}}(\mathbf{X})<\mathbf{X}<VaR_{\boldsymbol{q}}(\mathbf{X})\right]\\
&=\mathrm{E}\left[\mathbf{X}\mathbf{X}^{T}|VaR_{\boldsymbol{p}}(\mathbf{X})<\mathbf{X}<VaR_{\boldsymbol{q}}(\mathbf{X})\right]-\mathrm{MDTE}_{\boldsymbol{p},\boldsymbol{q}}(\mathbf{X})\mathrm{MDTE}_{\boldsymbol{p},\boldsymbol{q}}^{T}(\mathbf{X}).
\end{align*}
Using the transformation $\mathbf{Y}=\mathbf{\Sigma}^{-\frac{1}{2}}(\mathbf{Y}-\boldsymbol{\mu})$ and basic algebraic calculations, we obtain
\begin{align*}
&\mathrm{MDTCov}_{\boldsymbol{p},\boldsymbol{q}}(\mathbf{X})\\
&=\mathbf{\Sigma}^{\frac{1}{2}}\{\mathrm{E}[\boldsymbol{Y}\boldsymbol{Y}^{T}|\boldsymbol{\eta_{p}}<\mathbf{Y}<\boldsymbol{\eta_{q}}]-\mathrm{MDTE}_{(\boldsymbol{\eta_{p}},\boldsymbol{\eta_{q}})}(\mathbf{Y})\mathrm{MDTE}_{(\boldsymbol{\eta_{p}},\boldsymbol{\eta_{q}})}(\mathbf{Y})^{T}\}\mathbf{\Sigma}^{\frac{1}{2}},
\end{align*}
where $\boldsymbol{\eta_{v}}=\mathbf{\Sigma}^{-\frac{1}{2}}(VaR_{\boldsymbol{v}}(\mathbf{X})-\boldsymbol{\mu})$,~$\boldsymbol{v=p,q}$.\\
 Note that
\begin{align*}
&\mathrm{E}[Y_{i}Y_{j}|\boldsymbol{\eta_{p}}<\mathbf{Y}<\boldsymbol{\eta_{q}}]=\frac{1}{F_{\mathbf{Y}}(\boldsymbol{\eta_{p}},\boldsymbol{\eta_{q}})}\int_{\boldsymbol{\eta_{p}}}^{\boldsymbol{\eta_{q}}}y_{i}y_{j}c_{n}g_{n}\left(\frac{1}{2}\boldsymbol{y}^{T}\boldsymbol{y}\right)\mathrm{d}\boldsymbol{y}\\
&=\frac{c_{n}}{F_{\mathbf{Y}}(\boldsymbol{\eta_{p}},\boldsymbol{\eta_{q}})}\int_{\boldsymbol{\eta}_{\boldsymbol{p},-i}}^{\boldsymbol{\eta}_{\boldsymbol{q},-i}}y_{j}\int_{\eta_{\boldsymbol{p},i}}^{\eta_{\boldsymbol{q},i}}y_{i}g_{n}\left(\frac{1}{2}\boldsymbol{y}_{-i}^{T}\boldsymbol{y}_{-i}+\frac{1}{2}y_{i}^{2}\right)\mathrm{d}y_{i}\mathrm{d}\boldsymbol{y}_{-i},~for ~i\neq j,
\end{align*}
where $\boldsymbol{y}_{-i}=(y_{1},\cdots,y_{i-1},y_{i+1},\cdots,y_{n})^{T}$.\\
Since
\begin{align*}
&\int_{\eta_{\boldsymbol{p},i}}^{\eta_{\boldsymbol{q},i}}y_{i}g_{n}\left(\frac{1}{2}\boldsymbol{y}_{-i}^{T}\boldsymbol{y}_{-i}+\frac{1}{2}y_{i}^{2}\right)\mathrm{d}y_{i}\\
=&-\int_{\eta_{\boldsymbol{p},i}}^{\eta_{\boldsymbol{q},i}}\partial_{i}\overline{G}_{n}\left(\frac{1}{2}\boldsymbol{y}_{-i}^{T}\boldsymbol{y}_{-i}+\frac{1}{2}y_{i}^{2}\right)\\
=&\overline{G}_{n}\left(\frac{1}{2}\boldsymbol{y}_{-i}^{T}\boldsymbol{y}_{-i}+\frac{1}{2}\eta_{\boldsymbol{p},i}^{2}\right)-\overline{G}_{n}\left(\frac{1}{2}\boldsymbol{y}_{-i}^{T}\boldsymbol{y}_{-i}+\frac{1}{2}\eta_{\boldsymbol{q},i}^{2}\right),
\end{align*}
\begin{align*}
&\mathrm{E}[Y_{i}Y_{j}|\boldsymbol{\eta_{p}}<\mathbf{Y}<\boldsymbol{\eta_{q}}]\\
=&\frac{c_{n}}{F_{\mathbf{Y}}(\boldsymbol{\eta_{p}},\boldsymbol{\eta_{q}})}\int_{\boldsymbol{\eta}_{\boldsymbol{p},-i}}^{\boldsymbol{\eta}_{\boldsymbol{q},-i}}y_{j}\left[\overline{G}_{n}\left(\frac{1}{2}\boldsymbol{y}_{-i}^{T}\boldsymbol{y}_{-i}+\frac{1}{2}\eta_{\boldsymbol{p},i}^{2}\right)-\overline{G}_{n}\left(\frac{1}{2}\boldsymbol{y}_{-i}^{T}\boldsymbol{y}_{-i}+\frac{1}{2}\eta_{\boldsymbol{q},i}^{2}\right)\right]\mathrm{d}\boldsymbol{y}_{-i}\\
=&\frac{c_{n}}{F_{\mathbf{Y}}(\boldsymbol{\eta_{p}},\boldsymbol{\eta_{q}})}\int_{\boldsymbol{\eta}_{\boldsymbol{p},-ij}}^{\boldsymbol{\eta}_{\boldsymbol{q},-ij}}\int_{\eta_{\boldsymbol{p},j}}^{\eta_{\boldsymbol{q},j}}y_{j}\bigg[\overline{G}_{n}\left(\frac{1}{2}\boldsymbol{y}_{-ij}^{T}\boldsymbol{y}_{-ij}+\frac{1}{2}y_{j}^{2}+\frac{1}{2}\eta_{\boldsymbol{p},i}^{2}\right)\\
&-\overline{G}_{n}\left(\frac{1}{2}\boldsymbol{y}_{-ij}^{T}\boldsymbol{y}_{-ij}+\frac{1}{2}y_{j}^{2}+\frac{1}{2}\eta_{\boldsymbol{q},i}^{2}\right)\bigg]\mathrm{d}y_{j}\mathrm{d}\boldsymbol{y}_{-ij},~i\neq j,
\end{align*}
where $\boldsymbol{y}_{-ij}=(y_{1},\cdots,y_{i-1},y_{i+1},\cdots,y_{j-1},y_{j+1},\cdots,y_{n})^{T}$.\\
While
\begin{align*}
&\int_{\eta_{\boldsymbol{p},j}}^{\eta_{\boldsymbol{q},j}}y_{j}\overline{G}_{n}\left(\frac{1}{2}\boldsymbol{y}_{-ij}^{T}\boldsymbol{y}_{-ij}+\frac{1}{2}y_{j}^{2}+\frac{1}{2}\eta_{\boldsymbol{v},i}^{2}\right)\mathrm{d}y_{j}\\
&=-\int_{\eta_{\boldsymbol{p},j}}^{\eta_{\boldsymbol{q},j}}\partial_{j}\overline{\mathcal{G}}_{n}\left(\frac{1}{2}\boldsymbol{y}_{-ij}^{T}\boldsymbol{y}_{-ij}+\frac{1}{2}y_{j}^{2}+\frac{1}{2}\eta_{\boldsymbol{v},i}^{2}\right)\\
&=\overline{\mathcal{G}}_{n}\left(\frac{1}{2}\boldsymbol{y}_{-ij}^{T}\boldsymbol{y}_{-ij}+\frac{1}{2}\eta_{\boldsymbol{p},j}^{2}+\frac{1}{2}\eta_{\boldsymbol{v},i}^{2}\right)\\
&~~~-\overline{\mathcal{G}}_{n}\left(\frac{1}{2}\boldsymbol{y}_{-ij}^{T}\boldsymbol{y}_{-ij}+\frac{1}{2}\eta_{\boldsymbol{q},j}^{2}+\frac{1}{2}\eta_{\boldsymbol{v},i}^{2}\right)
,~\boldsymbol{v=p,q},
\end{align*}
thus
\begin{align*}
&\mathrm{E}[Y_{i}Y_{j}|\boldsymbol{\eta_{p}}<\mathbf{Y}<\boldsymbol{\eta_{q}}]=\frac{1}{F_{\mathbf{Y}}(\boldsymbol{\eta_{p}},\boldsymbol{\eta_{q}})}\\
& \boldsymbol{\cdot}\bigg\{\frac{c_{n}}{c_{n-2,\boldsymbol{p}i,\boldsymbol{p}j,ij}^{\ast\ast}}F_{\mathbf{Y}_{\boldsymbol{p}i,\boldsymbol{p}j,-ij}}(\boldsymbol{\eta}_{\boldsymbol{p},-ij},\boldsymbol{\eta}_{\boldsymbol{q},-ij})-\frac{c_{n}}{c_{n-2,\boldsymbol{p}i,\boldsymbol{q}j,ij}^{\ast\ast}}F_{\mathbf{Y}_{\boldsymbol{p}i,\boldsymbol{q}j,-ij}}(\boldsymbol{\eta}_{\boldsymbol{p},-ij},\boldsymbol{\eta}_{\boldsymbol{q},-ij})\\
&+\frac{c_{n}}{c_{n-2,\boldsymbol{q}i,\boldsymbol{q}j,ij}^{\ast\ast}}F_{\mathbf{Y}_{\boldsymbol{q}i,\boldsymbol{q}j,-ij}}(\boldsymbol{\eta}_{\boldsymbol{p},-ij},\boldsymbol{\eta}_{\boldsymbol{q},-ij}) -\frac{c_{n}}{c_{n-2,\boldsymbol{p}j,\boldsymbol{q}i,ij}^{\ast\ast}}F_{\mathbf{Y}_{\boldsymbol{p}j,\boldsymbol{q}i,-ij}}(\boldsymbol{\eta}_{\boldsymbol{p},-ij},\boldsymbol{\eta}_{\boldsymbol{q},-ij}) \bigg\},\\
&~i\neq j.
\end{align*}
In a similar manner, by using integration by parts, after some algebra we obtain
\begin{align*}
&\mathrm{E}[Y_{i}^{2}|\boldsymbol{\eta_{p}}<\mathbf{Z}<\boldsymbol{\eta_{q}}]=\frac{c_{n}}{F_{\mathbf{Y}}(\boldsymbol{\eta_{p}},\boldsymbol{\eta_{q}})}\int_{\boldsymbol{\eta_{p}}}^{\boldsymbol{\eta_{q}}}y_{i}^{2}g_{n}\left(\frac{1}{2}\boldsymbol{y}^{T}\boldsymbol{y}\right)\mathrm{d}\boldsymbol{y}\\
&=\frac{c_{n}}{F_{\mathbf{Y}}(\boldsymbol{\eta_{p}},\boldsymbol{\eta_{q}})}\int_{\boldsymbol{\eta}_{\boldsymbol{p},-i}}^{\boldsymbol{\eta}_{\boldsymbol{q},-i}}\int_{\eta_{\boldsymbol{p},i}}^{\eta_{\boldsymbol{q},i}}y_{i}^{2}g_{n}\left(\frac{1}{2}\boldsymbol{y}_{-i}^{T}\boldsymbol{y}_{-i}+\frac{1}{2}y_{i}^{2}\right)\mathrm{d}y_{i}\mathrm{d}\boldsymbol{y}_{-i}\\
&=\frac{c_{n}}{F_{\mathbf{Y}}(\boldsymbol{\eta_{p}},\boldsymbol{\eta_{q}})}\int_{\boldsymbol{\eta}_{\boldsymbol{p},-i}}^{\boldsymbol{\eta}_{\boldsymbol{q},-i}}\int_{\eta_{\boldsymbol{p},i}}^{\eta_{\boldsymbol{q},i}}-y_{i}\partial_{i}\overline{G}_{n}\left(\frac{1}{2}\boldsymbol{y}_{-i}^{T}\boldsymbol{y}_{-i}+\frac{1}{2}y_{i}^{2}\right)\mathrm{d}\boldsymbol{y}_{-i}\\
&=\frac{c_{n}}{F_{\mathbf{Y}}(\boldsymbol{\eta_{p}},\boldsymbol{\eta_{q}})}\int_{\boldsymbol{\eta}_{\boldsymbol{p},-i}}^{\boldsymbol{\eta}_{\boldsymbol{q},-i}}
\bigg[\eta_{\boldsymbol{p},i}\overline{G}_{n}\left(\frac{1}{2}\boldsymbol{y}_{-i}^{T}\boldsymbol{y}_{-i}+\frac{1}{2}\eta_{\boldsymbol{p},i}^{2}\right)-\eta_{\boldsymbol{q},i}\overline{G}_{n}\left(\frac{1}{2}\boldsymbol{y}_{-i}^{T}\boldsymbol{y}_{-i}+\frac{1}{2}\eta_{\boldsymbol{q},i}^{2}\right)\\
&~~~+\int_{\eta_{\boldsymbol{p},i}}^{\eta_{\boldsymbol{q},i}}\overline{G}_{n}\left(\frac{1}{2}\boldsymbol{y}_{-i}^{T}\boldsymbol{y}_{-i}+\frac{1}{2}y_{i}^{2}\right)\mathrm{d}y_{i}\bigg]
\mathrm{d}\boldsymbol{y}_{-i}\\
&=\frac{c_{n}}{F_{\mathbf{Y}}(\boldsymbol{\eta_{p}},\boldsymbol{\eta_{q}})}\bigg\{ \frac{\eta_{\boldsymbol{p},i}}{c_{n-1,\boldsymbol{p},i}}F_{\mathbf{Y}_{\boldsymbol{p},-i}}(\boldsymbol{\eta}_{\boldsymbol{p},-i},\boldsymbol{\eta}_{\boldsymbol{q},-i})-\frac{\eta_{\boldsymbol{q},i}}{c_{n-1,\boldsymbol{q},i}}F_{\mathbf{Y}_{\boldsymbol{q},-i}}(\boldsymbol{\eta}_{\boldsymbol{p},-i},\boldsymbol{\eta}_{\boldsymbol{q},-i})\\
&~~~+\frac{1}{c_{n}^{\ast}}F_{\mathbf{Y}^{\ast}}(\boldsymbol{\eta_{p}},\boldsymbol{\eta_{q}})\bigg\}.
\end{align*}
As for $\mathrm{MDTE}_{(\boldsymbol{\eta_{p}},\boldsymbol{\eta_{q}})}(\mathbf{Y})_{k}$, using Theorem 1 we immediately obtain (\ref{(37)}).
Therefore we obtain $(\ref{(36)})$, as required.

When $n=1$, we obtain DTV for univariate elliptical distribution:
\begin{align}\label{(38)}
&\mathrm{DTV}_{(p,q)}(X)=\frac{\sigma^{2}}{F_{Y}(\eta_{p},\eta_{q})}\left\{\lambda_{p,q}+\frac{c_{1}}{c_{1}^{\ast}}F_{Y^{\ast}}(\eta_{p},\eta_{q})-\delta_{p,q}^{2}\right\},
\end{align}
where
\begin{align*}
\lambda_{p,q}=c_{1}\left[\eta_{p}\overline{G}_{1}\left(\frac{1}{2}\eta_{p}^{2}\right)-\eta_{q}\overline{G}_{1}\left(\frac{1}{2}\eta_{q}^{2}\right)\right],
\end{align*}
$Y\sim E_{1}\left(0,~1,~g_{1}\right)$, $Y^{\ast}\sim E_{1}\left(0,~1,~\overline{G}_{1}\right)$, and $\delta_{p,q}$ is the same as in (\ref{(a33)}).\\
$\mathbf{Remark~6}$ We find that $\boldsymbol{x_{q}}=VaR_{\boldsymbol{q}}(\mathbf{X})\rightarrow\boldsymbol{+\infty}$ in Theorem $2$, one gets the formula of Theorem 2 in Landsman et al. (2018); When $x_{q}=VaR_{q}(X)\rightarrow +\infty$ in (\ref{(38)}), it is coincide with the result of (1.7) in  Furman and Landsman (2006), which is the generalization of Theorem 2 in Jiang et al. (2016).\\
\noindent$\mathbf{Corollary~6}$
  Let $\mathbf{X}\sim N_{n}(\boldsymbol{\mu},~\mathbf{\Sigma})$ $(n\geq2)$. Then
\begin{align}\label{(39)} &\mathrm{MDTCov}_{(\boldsymbol{p},\boldsymbol{q})}(\mathbf{X})=\mathbf{\Sigma}^{\frac{1}{2}}\mathbf{\Upsilon}_{\boldsymbol{p},\boldsymbol{q}}\mathbf{\Sigma}^{\frac{1}{2}},
\end{align}
where
\begin{align*}
&\mathbf{\Upsilon}_{\boldsymbol{p,q},i,j}= \frac{F_{\mathbf{Y},-ij}(\boldsymbol{\eta}_{\boldsymbol{p},-ij},\boldsymbol{\eta}_{\boldsymbol{q},-ij})}{ F_{\mathbf{Y}}(\boldsymbol{\eta_{p}},\boldsymbol{\eta_{q}})}\bigg\{\phi(\eta_{\boldsymbol{p},i})\phi(\eta_{\boldsymbol{p},j})-\phi(\eta_{\boldsymbol{p},i})\phi(\eta_{\boldsymbol{q},j})\\
&~~~+\phi(\eta_{\boldsymbol{q},i})\phi(\eta_{\boldsymbol{q},j})-\phi(\eta_{\boldsymbol{p},j})\phi(\eta_{\boldsymbol{q},i})\bigg\}-\mathrm{MDTE}_{(\boldsymbol{\eta_{p}},\boldsymbol{\eta_{q}})}(\mathbf{Y})_{i}\mathrm{MDTE}_{(\boldsymbol{\eta_{p}},\boldsymbol{\eta_{q}})}(\mathbf{Y})_{j}
,\\
&~~~i\neq j,~i=1,2,\cdots,n,~j=1,2,\cdots,n,
\end{align*}
and
\begin{align*}
\mathbf{\Upsilon}_{\boldsymbol{p,q},ii}
&=\frac{F_{\mathbf{Y}_{-i}}(\boldsymbol{\eta}_{\boldsymbol{p},-i},\boldsymbol{\eta}_{\boldsymbol{q},-i})}{F_{\mathbf{Y}}(\boldsymbol{\eta}_{\boldsymbol{p}},\boldsymbol{\eta}_{\boldsymbol{q}})}\left[\eta_{\boldsymbol{p},i}\phi(\eta_{\boldsymbol{p},i})-\eta_{\boldsymbol{q},i}\phi(\eta_{\boldsymbol{q},i})\right]\\
&~~~+1-[\mathrm{MDTE}_{(\boldsymbol{\eta_{p}},\boldsymbol{\eta_{q}})}(\mathbf{Y})_{i}]^{2},~i=1,2,\cdots,n,
\end{align*}
with
\begin{align}\label{(37)}
\mathrm{MDTE}_{(\boldsymbol{\eta_{p}},\boldsymbol{\eta_{q}})}(\mathbf{Y})_{k}
&=\frac{F_{\mathbf{Y}_{-k}}(\boldsymbol{\eta}_{\boldsymbol{p},-k},\boldsymbol{\eta}_{\boldsymbol{q},-k})}{F_{\mathbf{Y}}(\boldsymbol{\eta}_{\boldsymbol{p}},\boldsymbol{\eta}_{\boldsymbol{q}})}\left[\phi(\eta_{\boldsymbol{p},k})-\phi(\eta_{\boldsymbol{q},k})\right],
\end{align}
  $k=1,~2,\cdots,n,$ $\mathbf{Y}\sim N_{n}\left(\boldsymbol{0},~\boldsymbol{I_{n}}\right)$, $\mathbf{Y}_{-k}\sim N_{n-1}(\boldsymbol{0},~\mathbf{I_{n-1}})$, $\mathbf{Y}_{-ij}\sim N_{n-2}(\boldsymbol{0},~\mathbf{I_{n-2}})$.

When $n=1$, the DTV for univariate normal distribution is given:
\begin{align}\label{(40)}
&\mathrm{DTV}_{(p,q)}(X)=\frac{\sigma^{2}}{F_{Y}(\eta_{p},\eta_{q})}\left\{\lambda_{p,q}+F_{Y}(\eta_{p},\eta_{q})-\delta_{p,q}^{2}\right\},
\end{align}
where
\begin{align*}
\lambda_{p,q}=\eta_{p}\phi(\eta_{p})-\eta_{q}\phi(\eta_{q}),
\end{align*}
$Y\sim N_{1}\left(0,~1\right)$, and $\delta_{p,q}$ is the same as that in (\ref{(a25)}). \\
\noindent$\mathbf{Corollary~7}$
  Let $\mathbf{X}\sim St_{n}(\boldsymbol{\mu},~\mathbf{\Sigma},m)$ $(n\geq2)$. Then
\begin{align}\label{(41)} &\mathrm{MDTCov}_{(\boldsymbol{p},\boldsymbol{q})}(\mathbf{X})=\mathbf{\Sigma}^{\frac{1}{2}}\mathbf{\Upsilon}_{\boldsymbol{p},\boldsymbol{q}}\mathbf{\Sigma}^{\frac{1}{2}},
\end{align}
where
\begin{align*}
&\mathbf{\Upsilon}_{\boldsymbol{p,q},i,j}=\frac{1}{F_{\mathbf{Y}}(\boldsymbol{\eta_{p}},\boldsymbol{\eta_{q}})}\\
& \boldsymbol{\cdot}\bigg\{\frac{c_{n}}{c_{n-2,\boldsymbol{p}i,\boldsymbol{p}j,ij}^{\ast\ast}}F_{\mathbf{Y}_{\boldsymbol{p}i,\boldsymbol{p}j,-ij}}(\boldsymbol{\eta}_{\boldsymbol{p},-ij},\boldsymbol{\eta}_{\boldsymbol{q},-ij})-\frac{c_{n}}{c_{n-2,\boldsymbol{p}i,\boldsymbol{q}j,ij}^{\ast\ast}}F_{\mathbf{Y}_{\boldsymbol{p}i,\boldsymbol{q}j,-ij}}(\boldsymbol{\eta}_{\boldsymbol{p},-ij},\boldsymbol{\eta}_{\boldsymbol{q},-ij})\\
&+\frac{c_{n}}{c_{n-2,\boldsymbol{q}i,\boldsymbol{q}j,ij}^{\ast\ast}}F_{\mathbf{Y}_{\boldsymbol{q}i,\boldsymbol{q}j,-ij}}(\boldsymbol{\eta}_{\boldsymbol{p},-ij},\boldsymbol{\eta}_{\boldsymbol{q},-ij}) -\frac{c_{n}}{c_{n-2,\boldsymbol{p}j,\boldsymbol{q}i,ij}^{\ast\ast}}F_{\mathbf{Y}_{\boldsymbol{p}j,\boldsymbol{q}i,-ij}}(\boldsymbol{\eta}_{\boldsymbol{p},-ij},\boldsymbol{\eta}_{\boldsymbol{q},-ij}) \bigg\}\\
&-\mathrm{MDTE}_{(\boldsymbol{\eta_{p}},\boldsymbol{\eta_{q}})}(\mathbf{Y})_{i}\mathrm{MDTE}_{(\boldsymbol{\eta_{p}},\boldsymbol{\eta_{q}})}(\mathbf{Y})_{j}
,~i\neq j,~i=1,2,\cdots,n,~j=1,2,\cdots,n,
\end{align*}
and
\begin{align*}
&\mathbf{\Upsilon}_{\boldsymbol{p,q},ii}=\\ &\frac{1}{F_{\mathbf{Y}}(\boldsymbol{\eta}_{\boldsymbol{p}},\boldsymbol{\eta}_{\boldsymbol{q}})}\bigg\{\frac{c_{n}}{c_{n-1,\boldsymbol{p},i}^{\ast}}\eta_{\boldsymbol{p},i}F_{\mathbf{Y}_{\boldsymbol{p},-i}}(\boldsymbol{\eta}_{\boldsymbol{p},-i},\boldsymbol{\eta}_{\boldsymbol{q},-i}) -\frac{c_{n}}{c_{n-1,\boldsymbol{q},i}^{\ast}}\eta_{\boldsymbol{q},i}F_{\mathbf{Y}_{\boldsymbol{q},-i}}(\boldsymbol{\eta}_{\boldsymbol{p},-i},\boldsymbol{\eta}_{\boldsymbol{q},-i}) \\ &~~~+\frac{c_{n}}{c_{n}^{\ast}}F_{\mathbf{Y}^{\ast}}(\boldsymbol{\eta}_{\boldsymbol{p}},\boldsymbol{\eta}_{\boldsymbol{q}})\bigg\}-[\mathrm{MDTE}_{(\boldsymbol{\eta_{p}},\boldsymbol{\eta_{q}})}(\mathbf{Y})_{i}]^{2},~i=1,2,\cdots,n,
\end{align*}
with
\begin{align}\label{(42)}
\nonumber&\mathrm{MDTE}_{(\boldsymbol{\eta_{p}},\boldsymbol{\eta_{q}})}(\mathbf{Y})_{k}\\
&=\frac{1}{F_{\mathbf{Y}}(\boldsymbol{\eta}_{\boldsymbol{p}},\boldsymbol{\eta}_{\boldsymbol{q}})}\bigg\{\frac{c_{n}}{c_{n-1,\boldsymbol{p},k}^{\ast}}F_{\mathbf{Y}_{\boldsymbol{p},-k}}(\boldsymbol{\eta}_{\boldsymbol{p},-k},\boldsymbol{\eta}_{\boldsymbol{q},-k})-\frac{c_{n}}{c_{n-1,\boldsymbol{q},k}^{\ast}}F_{\mathbf{Y}_{\boldsymbol{q},-k}}(\boldsymbol{\eta}_{\boldsymbol{p},-k},\boldsymbol{\eta}_{\boldsymbol{q},-k})\bigg\},
\end{align}
  $k=1,~2,\cdots,n,$ $\mathbf{Y}\sim St_{n}\left(\boldsymbol{0},~\boldsymbol{I_{n}},~m\right)$, $\mathbf{Y}^{\ast}\sim E_{n}\left(\boldsymbol{0},~\boldsymbol{I_{n}},~\overline{G}_{n}\right)$,\\
  $\mathbf{Y}_{\boldsymbol{u}r,\boldsymbol{v}l,-ij}\sim St_{n-2}\left(\boldsymbol{0},~\Delta_{n-2,\boldsymbol{u}r,\boldsymbol{v}l,ij},~m-2\right)$, $\mathbf{Y}_{\boldsymbol{v},-k}\sim St_{n-1}\left(\boldsymbol{0},~\Delta_{\boldsymbol{v},k},~m-1\right),$
 $$\Delta_{n-2,\boldsymbol{u}r,\boldsymbol{v}l,ij}=\left(\frac{m+\eta_{\boldsymbol{u}r}^{2}+\eta_{\boldsymbol{v}l}^{2}}{m-2}\right)\mathbf{I}_{n-2},$$ $$c_{n-2,\boldsymbol{u}r,\boldsymbol{v}l,ij}^{\ast\ast}=\frac{\Gamma\left(\frac{m+n-4}{2}\right)(m+n-2)(m+n-4)}{\Gamma\left(\frac{m-2}{2}\right)\pi^{\frac{n-2}{2}}m^{\frac{m+n}{2}}}\left(m+\eta_{\boldsymbol{u}r}^{2}+\eta_{\boldsymbol{v}l}^{2}\right)^{\frac{m+n-4}{2}},$$ $\boldsymbol{u}=\boldsymbol{p},\boldsymbol{q};\boldsymbol{v}=\boldsymbol{p},\boldsymbol{q};r=i,j;l=i,j,$
 and $\Delta_{\boldsymbol{v},k}$ is the same as that in Corollary 2. In addition, $c_{n}$, $c_{n}^{\ast}$ and $\overline{G}_{n}$ are in (\ref{(d18)}), (\ref{(e18)}) and (\ref{(b18)}), respectively.

 We can further simplify
 $$\frac{c_{n}}{c_{n-2,\boldsymbol{u}r,\boldsymbol{v}l,ij}^{\ast\ast}}=\frac{(m+n)m^{\frac{m-2}{2}}}{2(m+n-4)\pi}\left(m+\eta_{\boldsymbol{u}r}^{2}+\eta_{\boldsymbol{v}l}^{2}\right)^{-(m+n-4)/2}$$ and
 $$
\frac{c_{n}}{c_{n}^{\ast}} =\frac{p\Gamma(\frac{m+n}{2})B(\frac{n}{2},~\frac{m-2}{2})}{(m+n-2)\Gamma(\frac{m}{2})\Gamma(\frac{n}{2})}=\frac{m}{m-2},~\mathrm{if}~m>2.
$$

When $n=1$, we obtain the DTV for univariate student-$t$ distribution:
\begin{align}\label{(43)}
&\mathrm{DTV}_{(p,q)}(X)=\frac{\sigma^{2}}{F_{Y}(\eta_{p},\eta_{q})}\left\{\lambda_{p,q}+\frac{m}{m-2}F_{Y^{\ast}}(\eta_{p},\eta_{q})-\delta_{p,q}^{2}\right\},~m>2,
\end{align}
where
\begin{align*}
\lambda_{p,q}=\frac{m\Gamma((m+1)/2)}{(m-1)\Gamma(m/2)\sqrt{m\pi}}\left[\eta_{p}\left(1+\frac{\eta_{p}^{2}}{m}\right)^{-(m-1)/2}-\eta_{q}\left(1+\frac{\eta_{q}^{2}}{m}\right)^{-(m-1)/2}\right],
\end{align*}
$Y\sim St_{1}\left(0,~1,~m\right)$, $Y^{\ast}\sim E_{1}\left(0,~1,~\overline{G}_{1}\right)$, and $\overline{G}_{1}$ and $\delta_{p,q}$ are in (\ref{(b18)}) and (\ref{(a27)}), respectively.\\
\noindent$\mathbf{Corollary~8}$
  Let $\mathbf{X}\sim Lo_{n}(\boldsymbol{\mu},~\mathbf{\Sigma})$ $(n\geq2)$. Then
\begin{align}\label{(44)} &\mathrm{MDTCov}_{(\boldsymbol{p},\boldsymbol{q})}(\mathbf{X})=\mathbf{\Sigma}^{\frac{1}{2}}\mathbf{\Upsilon}_{\boldsymbol{p},\boldsymbol{q}}\mathbf{\Sigma}^{\frac{1}{2}},
\end{align}
where
 \begin{align*}
&\mathbf{\Upsilon}_{\boldsymbol{p,q},i,j}= \frac{1}{F_{\mathbf{Y}}(\boldsymbol{\eta_{p}},\boldsymbol{\eta_{q}})}\\
&\boldsymbol{\cdot}\bigg\{\frac{c_{n}}{c_{n-2,\boldsymbol{p}i,\boldsymbol{p}j,ij}^{\ast\ast}}F_{\mathbf{Y}_{\boldsymbol{p}i,\boldsymbol{p}j,-ij}}(\boldsymbol{\eta}_{\boldsymbol{p},-ij},\boldsymbol{\eta}_{\boldsymbol{q},-ij})-\frac{c_{n}}{c_{n-2,\boldsymbol{p}i,\boldsymbol{q}j,ij}^{\ast\ast}}F_{\mathbf{Y}_{\boldsymbol{p}i,\boldsymbol{q}j,-ij}}(\boldsymbol{\eta}_{\boldsymbol{p},-ij},\boldsymbol{\eta}_{\boldsymbol{q},-ij})\\
&+\frac{c_{n}}{c_{n-2,\boldsymbol{q}i,\boldsymbol{q}j,ij}^{\ast\ast}}F_{\mathbf{Y}_{\boldsymbol{q}i,\boldsymbol{q}j,-ij}}(\boldsymbol{\eta}_{\boldsymbol{p},-ij},\boldsymbol{\eta}_{\boldsymbol{q},-ij}) -\frac{c_{n}}{c_{n-2,\boldsymbol{p}j,\boldsymbol{q}i,ij}^{\ast\ast}}F_{\mathbf{Y}_{\boldsymbol{p}j,\boldsymbol{q}i,-ij}}(\boldsymbol{\eta}_{\boldsymbol{p},-ij},\boldsymbol{\eta}_{\boldsymbol{q},-ij}) \bigg\}\\
&-\mathrm{MDTE}_{(\boldsymbol{\eta_{p}},\boldsymbol{\eta_{q}})}(\mathbf{Y})_{i}\mathrm{MDTE}_{(\boldsymbol{\eta_{p}},\boldsymbol{\eta_{q}})}(\mathbf{Y})_{j}
,~i\neq j,~i=1,2,\cdots,n,~j=1,2,\cdots,n,
\end{align*}
and
\begin{align*}
&\mathbf{\Upsilon}_{\boldsymbol{p,q},ii}=\\ &\frac{1}{F_{\mathbf{Y}}(\boldsymbol{\eta}_{\boldsymbol{p}},\boldsymbol{\eta}_{\boldsymbol{q}})}\bigg\{\frac{c_{n}}{c_{n-1,\boldsymbol{p},i}^{\ast}}\eta_{\boldsymbol{p},i}F_{\mathbf{Y}_{\boldsymbol{p},-i}}(\boldsymbol{\eta}_{\boldsymbol{p},-i},\boldsymbol{\eta}_{\boldsymbol{q},-i}) -\frac{c_{n}}{c_{n-1,\boldsymbol{q},i}^{\ast}}\eta_{\boldsymbol{q},i}F_{\mathbf{Y}_{\boldsymbol{q},-i}}(\boldsymbol{\eta}_{\boldsymbol{p},-i},\boldsymbol{\eta}_{\boldsymbol{q},-i}) \\ &~~~+\frac{\Psi_{1}^{\ast}(-1,\frac{n}{2},1)}{\Psi_{2}^{\ast}(-1,\frac{n}{2},1)}F_{\mathbf{Y}^{\ast}}(\boldsymbol{\eta}_{\boldsymbol{p}},\boldsymbol{\eta}_{\boldsymbol{q}})\bigg\}-[\mathrm{MDTE}_{(\boldsymbol{\eta_{p}},\boldsymbol{\eta_{q}})}(\mathbf{Y})_{i}]^{2},~i=1,2,\cdots,n,
\end{align*}
with
\begin{align}\label{(45)}
\nonumber&\mathrm{MDTE}_{(\boldsymbol{\eta_{p}},\boldsymbol{\eta_{q}})}(\mathbf{Y})_{k}\\
&=\frac{1}{F_{\mathbf{Y}}(\boldsymbol{\eta}_{\boldsymbol{p}},\boldsymbol{\eta}_{\boldsymbol{q}})}\bigg\{\frac{c_{n}}{c_{n-1,\boldsymbol{p},k}^{\ast}}F_{\mathbf{Y}_{\boldsymbol{p},-k}}(\boldsymbol{\eta}_{\boldsymbol{p},-k},\boldsymbol{\eta}_{\boldsymbol{q},-k})-\frac{c_{n}}{c_{n-1,\boldsymbol{q},k}^{\ast}}F_{\mathbf{Y}_{\boldsymbol{q},-k}}(\boldsymbol{\eta}_{\boldsymbol{p},-k},\boldsymbol{\eta}_{\boldsymbol{q},-k})\bigg\},
\end{align}
  $k=1,~2,\cdots,n,$ $\mathbf{Y}\sim Lo_{n}\left(\boldsymbol{0},~\boldsymbol{I_{n}}\right)$, $\mathbf{Y}^{\ast}\sim E_{n}\left(\boldsymbol{0},~\boldsymbol{I_{n}},~\overline{G}_{n}\right)$,\\ the pdf of $\mathbf{Y}_{\boldsymbol{u}r,\boldsymbol{v}l,-ij}$:
  $$f_{\mathbf{Y}_{\boldsymbol{u}r,\boldsymbol{v}l,-ij}}(\boldsymbol{t})=c_{n-2,\boldsymbol{u}r,\boldsymbol{v}l,ij}^{\ast\ast}\ln\left[1+\exp\left(-\frac{1}{2}\boldsymbol{t}^{T}\boldsymbol{t}-\frac{1}{2}\eta_{\boldsymbol{u},r}^{2}-\frac{1}{2}\eta_{\boldsymbol{v},l}^{2}\right)\right],~\boldsymbol{t}\in \mathbb{R}^{n-2},$$
  the normalizing constant
  \begin{align*}
&c_{n-2,\boldsymbol{u}r,\boldsymbol{v}l,ij}^{\ast\ast}\\
&=\frac{\Gamma((n-2)/2)}{(2\pi)^{(n-2)/2}}\left\{\int_{0}^{\infty}t^{(n-4)/2}\ln\left[1+\exp\left(-\frac{1}{2}\eta_{\boldsymbol{u},r}^{2}-\frac{1}{2}\eta_{\boldsymbol{v},l}^{2}\right)\exp(-t)\right]\mathrm{d}t\right\}^{-1},
\end{align*}
$\boldsymbol{u}=\boldsymbol{p},\boldsymbol{q};\boldsymbol{v}=\boldsymbol{p},\boldsymbol{q};r=i,j;l=i,j,$
and $\mathbf{Y}_{\boldsymbol{v},-k}$ and $c_{n-1,\boldsymbol{p},k}^{\ast}$ are the same as those in Corollary 3. In addition, $c_{n}$ and $\overline{G}_{n}$ are in (\ref{(d19)}) and (\ref{(b19)}), respectively.

We can further simplify
$$\frac{c_{n}}{c_{n-2,\boldsymbol{u}r,\boldsymbol{v}l,ij}^{\ast\ast}}=\frac{\int_{0}^{\infty}t^{(n-4)/2}\ln\left[1+\exp\left(-\frac{1}{2}\eta_{\boldsymbol{u},r}^{2}-\frac{1}{2}\eta_{\boldsymbol{v},l}^{2}\right)\exp(-t)\right]\mathrm{d}t}{2\Gamma((n-1)/2)\pi\Psi_{2}^{\ast}(-1,\frac{n}{2},1)}.$$

 When $n=1$, the DTV for univariate logistic distribution is given:
\begin{align}\label{(46)}
&\mathrm{DTV}_{(p,q)}(X)=\frac{\sigma^{2}}{F_{Y}(\eta_{p},\eta_{q})}\left\{\lambda_{p,q}+\frac{\Psi_{1}^{\ast}\left(-1,\frac{1}{2},1\right)}{\Psi_{2}^{\ast}\left(-1,\frac{1}{2},1\right)}F_{Y^{\ast}}(\eta_{p},\eta_{q})-\delta_{p,q}^{2}\right\},
\end{align}
where
\begin{align*}
\lambda_{p,q}=\frac{1}{\Psi_{2}^{\ast}\left(-1,\frac{1}{2},1\right)}\left[\eta_{p}\frac{\phi(\eta_{p})}{1+\sqrt{2\pi}\phi(\eta_{p})}-\eta_{q}\frac{\phi(\eta_{q})}{1+\sqrt{2\pi}\phi(\eta_{q})}\right],
\end{align*}
$Y\sim Lo_{1}\left(0,~1\right)$, $Y^{\ast}\sim E_{1}\left(0,~1,~\overline{G}_{1}\right)$, and $\overline{G}_{1}$ and $\delta_{p,q}$ are in (\ref{(b19)}) and (\ref{(a29)}), respectively.\\
 \noindent$\mathbf{Corollary~9}$
  Let $\mathbf{X}\sim La_{n}(\boldsymbol{\mu},~\mathbf{\Sigma})$ $(n\geq2)$. Then
\begin{align}\label{(47)} &\mathrm{MDTCov}_{(\boldsymbol{p},\boldsymbol{q})}(\mathbf{X})=\mathbf{\Sigma}^{\frac{1}{2}}\mathbf{\Upsilon}_{\boldsymbol{p},\boldsymbol{q}}\mathbf{\Sigma}^{\frac{1}{2}},
\end{align}
where
\begin{align*}
&\mathbf{\Upsilon}_{\boldsymbol{p,q},i,j}=\frac{1}{F_{\mathbf{Y}}(\boldsymbol{\eta_{p}},\boldsymbol{\eta_{q}})}\\
& \boldsymbol{\cdot}\bigg\{\frac{c_{n}}{c_{n-2,\boldsymbol{p}i,\boldsymbol{p}j,ij}^{\ast\ast}}F_{\mathbf{Y}_{\boldsymbol{p}i,\boldsymbol{p}j,-ij}}(\boldsymbol{\eta}_{\boldsymbol{p},-ij},\boldsymbol{\eta}_{\boldsymbol{q},-ij})-\frac{c_{n}}{c_{n-2,\boldsymbol{p}i,\boldsymbol{q}j,ij}^{\ast\ast}}F_{\mathbf{Y}_{\boldsymbol{p}i,\boldsymbol{q}j,-ij}}(\boldsymbol{\eta}_{\boldsymbol{p},-ij},\boldsymbol{\eta}_{\boldsymbol{q},-ij})\\
&+\frac{c_{n}}{c_{n-2,\boldsymbol{q}i,\boldsymbol{q}j,ij}^{\ast\ast}}F_{\mathbf{Y}_{\boldsymbol{q}i,\boldsymbol{q}j,-ij}}(\boldsymbol{\eta}_{\boldsymbol{p},-ij},\boldsymbol{\eta}_{\boldsymbol{q},-ij}) -\frac{c_{n}}{c_{n-2,\boldsymbol{p}j,\boldsymbol{q}i,ij}^{\ast\ast}}F_{\mathbf{Y}_{\boldsymbol{p}j,\boldsymbol{q}i,-ij}}(\boldsymbol{\eta}_{\boldsymbol{p},-ij},\boldsymbol{\eta}_{\boldsymbol{q},-ij}) \bigg\}\\
&-\mathrm{MDTE}_{(\boldsymbol{\eta_{p}},\boldsymbol{\eta_{q}})}(\mathbf{Y})_{i}\mathrm{MDTE}_{(\boldsymbol{\eta_{p}},\boldsymbol{\eta_{q}})}(\mathbf{Y})_{j}
,~i\neq j,~i=1,2,\cdots,n,~j=1,2,\cdots,n,
\end{align*}
and
\begin{align*}
&\mathbf{\Upsilon}_{\boldsymbol{p,q},ii}=\\ &\frac{1}{F_{\mathbf{Y}}(\boldsymbol{\eta}_{\boldsymbol{p}},\boldsymbol{\eta}_{\boldsymbol{q}})}\bigg\{\frac{c_{n}}{c_{n-1,\boldsymbol{p},i}^{\ast}}\eta_{\boldsymbol{p},i}F_{\mathbf{Y}_{\boldsymbol{p},-i}}(\boldsymbol{\eta}_{\boldsymbol{p},-i},\boldsymbol{\eta}_{\boldsymbol{q},-i}) -\frac{c_{n}}{c_{n-1,\boldsymbol{q},i}^{\ast}}\eta_{\boldsymbol{q},i}F_{\mathbf{Y}_{\boldsymbol{q},-i}}(\boldsymbol{\eta}_{\boldsymbol{p},-i},\boldsymbol{\eta}_{\boldsymbol{q},-i}) \\ &~~~+(n+1)F_{\mathbf{Y}^{\ast}}(\boldsymbol{\eta}_{\boldsymbol{p}},\boldsymbol{\eta}_{\boldsymbol{q}})\bigg\}-[\mathrm{MDTE}_{(\boldsymbol{\eta_{p}},\boldsymbol{\eta_{q}})}(\mathbf{Y})_{i}]^{2},~i=1,2,\cdots,n,
\end{align*}
with
\begin{align}\label{(48)}
\nonumber&\mathrm{MDTE}_{(\boldsymbol{\eta_{p}},\boldsymbol{\eta_{q}})}(\mathbf{Y})_{k}\\
&=\frac{1}{F_{\mathbf{Y}}(\boldsymbol{\eta}_{\boldsymbol{p}},\boldsymbol{\eta}_{\boldsymbol{q}})}\bigg\{\frac{c_{n}}{c_{n-1,\boldsymbol{p},k}^{\ast}}F_{\mathbf{Y}_{\boldsymbol{p},-k}}(\boldsymbol{\eta}_{\boldsymbol{p},-k},\boldsymbol{\eta}_{\boldsymbol{q},-k})-\frac{c_{n}}{c_{n-1,\boldsymbol{q},k}^{\ast}}F_{\mathbf{Y}_{\boldsymbol{q},-k}}(\boldsymbol{\eta}_{\boldsymbol{p},-k},\boldsymbol{\eta}_{\boldsymbol{q},-k})\bigg\},
\end{align}
  $k=1,~2,\cdots,n,$ $\mathbf{Y}\sim La_{n}\left(\boldsymbol{0},~\boldsymbol{I_{n}}\right)$, $\mathbf{Y}^{\ast}\sim E_{n}\left(\boldsymbol{0},~\boldsymbol{I_{n}},~\overline{G}_{n}\right)$,\\
  the pdf of $\mathbf{Y}_{\boldsymbol{u}r,\boldsymbol{v}l,-ij}$:
  \begin{align*}
f_{\mathbf{Y}_{\boldsymbol{u}r,\boldsymbol{v}l,-ij}}(\boldsymbol{t})=&c_{n-2,\boldsymbol{u}r,\boldsymbol{v}l,ij}^{\ast\ast}\left[\left(1+\frac{3}{\sqrt{2}}\right)\left(\boldsymbol{t}^{T}\boldsymbol{t}+\eta_{\boldsymbol{u},r}^{2}+\eta_{\boldsymbol{v},l}^{2}\right)+3\right]\\
&\cdot\exp(-\sqrt{\boldsymbol{t}^{T}\boldsymbol{t}+\eta_{\boldsymbol{u},r}^{2}+\eta_{\boldsymbol{v},l}^{2}}),~\boldsymbol{t}\in \mathbb{R}^{n-2},
\end{align*}
normalizing constant
 \begin{align*}
c_{n-2,\boldsymbol{u}r,\boldsymbol{v}l,ij}^{\ast\ast}=&\frac{\Gamma\left(\frac{n-2}{2}\right)}{(2\pi)^{(n-2)/2}}\bigg\{\int_{0}^{\infty}t^{\frac{n-4}{2}}\left[(2+3\sqrt{2})\left(t+\frac{\eta_{\boldsymbol{u},r}^{2}}{2}+\frac{\eta_{\boldsymbol{v},l}^{2}}{2}\right)+3\right]\\
&\cdot\exp\left(-\sqrt{2t+\eta_{\boldsymbol{u},r}^{2}+\eta_{\boldsymbol{v},l}^{2}}\right)\mathrm{d}t\bigg\}^{-1}.
\end{align*}
  $\boldsymbol{u}=\boldsymbol{p},\boldsymbol{q};\boldsymbol{v}=\boldsymbol{p},\boldsymbol{q};r=i,j;l=i,j,$
and $c_{n-1,\boldsymbol{v},k}^{\ast}$ and $\mathbf{Y}_{\boldsymbol{v},-k}$ are the same as those in Corollary 4. In addition, $c_{n}$ and $\overline{G}_{n}$ are in (\ref{(d20)}) and (\ref{(b20)}), respectively.

We can further simplify
\begin{align*}
\frac{c_{n}}{c_{n-2,\boldsymbol{u}r,\boldsymbol{v}l,ij}^{\ast\ast}}=&\frac{2^{(n-6)/2}}{\Gamma(n-1)\pi}\bigg\{\int_{0}^{\infty}t^{\frac{n-4}{2}}\left[(2+3\sqrt{2})\left(t+\frac{\eta_{\boldsymbol{u},r}^{2}}{2}+\frac{\eta_{\boldsymbol{v},l}^{2}}{2}\right)+3\right]\\
&\boldsymbol{\cdot}\exp\left(-\sqrt{2t+\eta_{\boldsymbol{u},r}^{2}+\eta_{\boldsymbol{v},l}^{2}}\right)\mathrm{d}t\bigg\}.
\end{align*}

When $n=1$, we obtain the DTV for univariate Laplace distribution:
\begin{align}\label{(49)}
&\mathrm{DTV}_{(p,q)}(X)=\frac{\sigma^{2}}{F_{Y}(\eta_{p},\eta_{q})}\left\{\lambda_{p,q}+2F_{Y^{\ast}}(\eta_{p},\eta_{q})-\delta_{p,q}^{2}\right\},
\end{align}
where
\begin{align*}
\lambda_{p,q}=\frac{1}{2}\left[\eta_{p}(1+|\eta_{p}|)\exp\left(-|\eta_{p}|\right)-\eta_{q}(1+|\eta_{q}|)\exp\left(-|\eta_{q}|\right)\right],
\end{align*}
$Y\sim La_{1}\left(0,~1\right)$, $Y^{\ast}\sim E_{1}\left(0,~1,~\overline{G}_{1}\right)$, and $\overline{G}_{1}$ and $\delta_{p,q}$ are in (\ref{(b20)}) and (\ref{(a31)}), respectively.\\
\noindent$\mathbf{Corollary~10}$
  Let $\mathbf{X}\sim PVII_{n}(\boldsymbol{\mu},~\mathbf{\Sigma},t)$ $(n\geq2)$. Then
\begin{align}\label{(50)} &\mathrm{MDTCov}_{(\boldsymbol{p},\boldsymbol{q})}(\mathbf{X})=\mathbf{\Sigma}^{\frac{1}{2}}\mathbf{\Upsilon}_{\boldsymbol{p},\boldsymbol{q}}\mathbf{\Sigma}^{\frac{1}{2}},
\end{align}
where
\begin{align*}
&\mathbf{\Upsilon}_{\boldsymbol{p,q},i,j}=\frac{1}{F_{\mathbf{Y}}(\boldsymbol{\eta_{p}},\boldsymbol{\eta_{q}})}\\
&\boldsymbol{\cdot} \bigg\{\frac{c_{n}}{c_{n-2,\boldsymbol{p}i,\boldsymbol{p}j,ij}^{\ast\ast}}F_{\mathbf{Y}_{\boldsymbol{p}i,\boldsymbol{p}j,-ij}}(\boldsymbol{\eta}_{\boldsymbol{p},-ij},\boldsymbol{\eta}_{\boldsymbol{q},-ij})-\frac{c_{n}}{c_{n-2,\boldsymbol{p}i,\boldsymbol{q}j,ij}^{\ast\ast}}F_{\mathbf{Y}_{\boldsymbol{p}i,\boldsymbol{q}j,-ij}}(\boldsymbol{\eta}_{\boldsymbol{p},-ij},\boldsymbol{\eta}_{\boldsymbol{q},-ij})\\
&+\frac{c_{n}}{c_{n-2,\boldsymbol{q}i,\boldsymbol{q}j,ij}^{\ast\ast}}F_{\mathbf{Y}_{\boldsymbol{q}i,\boldsymbol{q}j,-ij}}(\boldsymbol{\eta}_{\boldsymbol{p},-ij},\boldsymbol{\eta}_{\boldsymbol{q},-ij}) -\frac{c_{n}}{c_{n-2,\boldsymbol{p}j,\boldsymbol{q}i,ij}^{\ast\ast}}F_{\mathbf{Y}_{\boldsymbol{p}j,\boldsymbol{q}i,-ij}}(\boldsymbol{\eta}_{\boldsymbol{p},-ij},\boldsymbol{\eta}_{\boldsymbol{q},-ij}) \bigg\}\\
&-\mathrm{MDTE}_{(\boldsymbol{\eta_{p}},\boldsymbol{\eta_{q}})}(\mathbf{Y})_{i}\mathrm{MDTE}_{(\boldsymbol{\eta_{p}},\boldsymbol{\eta_{q}})}(\mathbf{Y})_{j}
,~i\neq j,~i=1,2,\cdots,n,~j=1,2,\cdots,n,
\end{align*}
and
\begin{align*}
&\mathbf{\Upsilon}_{\boldsymbol{p,q},ii}=\\ &\frac{1}{F_{\mathbf{Y}}(\boldsymbol{\eta}_{\boldsymbol{p}},\boldsymbol{\eta}_{\boldsymbol{q}})}\bigg\{\frac{c_{n}}{c_{n-1,\boldsymbol{p},i}^{\ast}}\eta_{\boldsymbol{p},i}F_{\mathbf{Y}_{\boldsymbol{p},-i}}(\boldsymbol{\eta}_{\boldsymbol{p},-i},\boldsymbol{\eta}_{\boldsymbol{q},-i}) -\frac{c_{n}}{c_{n-1,\boldsymbol{q},i}^{\ast}}\eta_{\boldsymbol{q},i}F_{\mathbf{Y}_{\boldsymbol{q},-i}}(\boldsymbol{\eta}_{\boldsymbol{p},-i},\boldsymbol{\eta}_{\boldsymbol{q},-i}) \\ &~~~+\frac{c_{n}}{c_{n}^{\ast}}F_{\mathbf{Y}^{\ast}}(\boldsymbol{\eta}_{\boldsymbol{p}},\boldsymbol{\eta}_{\boldsymbol{q}})\bigg\}-[\mathrm{MDTE}_{(\boldsymbol{\eta_{p}},\boldsymbol{\eta_{q}})}(\mathbf{Y})_{i}]^{2},~i=1,2,\cdots,n,
\end{align*}
with
\begin{align}\label{(51)}
\nonumber&\mathrm{MDTE}_{(\boldsymbol{\eta_{p}},\boldsymbol{\eta_{q}})}(\mathbf{Y})_{k}\\
&=\frac{1}{F_{\mathbf{Y}}(\boldsymbol{\eta}_{\boldsymbol{p}},\boldsymbol{\eta}_{\boldsymbol{q}})}\bigg\{\frac{c_{n}}{c_{n-1,\boldsymbol{p},k}^{\ast}}F_{\mathbf{Y}_{\boldsymbol{p},-k}}(\boldsymbol{\eta}_{\boldsymbol{p},-k},\boldsymbol{\eta}_{\boldsymbol{q},-k})-\frac{c_{n}}{c_{n-1,\boldsymbol{q},k}^{\ast}}F_{\mathbf{Y}_{\boldsymbol{q},-k}}(\boldsymbol{\eta}_{\boldsymbol{p},-k},\boldsymbol{\eta}_{\boldsymbol{q},-k})\bigg\},
\end{align}
  $k=1,~2,\cdots,n,$ $\mathbf{Y}\sim PVII_{n}\left(\boldsymbol{0},~\boldsymbol{I_{n}},~m\right)$, $\mathbf{Y}^{\ast}\sim E_{n}\left(\boldsymbol{0},~\boldsymbol{I_{n}},~\overline{G}_{n}\right)$,
  $$\mathbf{Y}_{\boldsymbol{u}r,\boldsymbol{v}l,-ij}\sim PVII_{n-2}\left(\boldsymbol{0},~\Lambda_{n-2,\boldsymbol{u}r,\boldsymbol{v}l,ij},~t-2\right),$$ $\mathbf{Y}_{\boldsymbol{v},-k}\sim PVII_{n-1}\left(\boldsymbol{0},~\Lambda_{\boldsymbol{v},k},~t-1\right),$
 $\Lambda_{n-2,\boldsymbol{u}r,\boldsymbol{v}l,ij}=\left(1+\eta_{\boldsymbol{u}r}^{2}+\eta_{\boldsymbol{v}l}^{2}\right)\mathbf{I}_{n-2},$ $$c_{n-2,\boldsymbol{u}r,\boldsymbol{v}l,ij}^{\ast\ast}=\frac{4\Gamma\left(t\right)}{\Gamma\left(t-\frac{n+2}{2}\right)\pi^{\frac{n-2}{2}}}\left(1+\eta_{\boldsymbol{u}r}^{2}+\eta_{\boldsymbol{v}l}^{2}\right)^{t-2},$$ $\boldsymbol{u}=\boldsymbol{p},\boldsymbol{q};\boldsymbol{v}=\boldsymbol{p},\boldsymbol{q};r=i,j;l=i,j,$
 and $\mathbf{Y}_{\boldsymbol{v},-k}$, $c_{n-1,\boldsymbol{v},k}^{\ast}$ and $\Lambda_{\boldsymbol{v},k}$ are the same as those in Corollary 5. In addition, $c_{n}$, $c_{n}^{\ast}$ and $\overline{G}_{n}$ are in (\ref{(d21)}), (\ref{(e21)}) and (\ref{(b21)}), respectively.

  We can further simplify
 $$\frac{c_{n}}{c_{n-2,\boldsymbol{u}r,\boldsymbol{v}l,ij}^{\ast\ast}}=\frac{1}{4(t-(n+2)/2)\pi}\left(1+\eta_{\boldsymbol{u}r}^{2}+\eta_{\boldsymbol{v}l}^{2}\right)^{-(t-2)}$$ and
 $$
\frac{c_{n}}{c_{n}^{\ast}} =\frac{1}{2t-n-2},~\mathrm{if}~t>n/2+1.
$$

 When $n=1$, the DTV for univariate Pearson type VII distribution is given:
\begin{align}\label{(52)}
&\mathrm{DTV}_{(p,q)}(X)=\frac{\sigma^{2}}{F_{Y}(\eta_{p},\eta_{q})}\left\{\lambda_{p,q}+\frac{1}{2t-3}F_{Y^{\ast}}(\eta_{p},\eta_{q})-\delta_{p,q}^{2}\right\},~t>\frac{3}{2},
\end{align}
where
\begin{align*}
\lambda_{p,q}=\frac{\Gamma(t)}{2(t-1)\Gamma(t-1/2)\sqrt{\pi}}\left[\eta_{p}\left(1+\eta_{p}^{2}\right)^{-(t-1)}-\eta_{q}\left(1+\eta_{q}^{2}\right)^{-(t-1)}\right],
\end{align*}
$Y\sim PVII_{1}\left(0,~1,~t\right)$, $Y^{\ast}\sim E_{1}\left(0,~1,~\overline{G}_{1}\right)$, and $\overline{G}_{1}$ and $\delta_{p,q}$ are in (\ref{(b21)}) and (\ref{(a27)}), respectively.
 \section{Numerical illustration}
 We consider DTE, DTV, MDTE and MDTCov risk measures for the normal (N) distributions.\\
  Let $\mathbf{X}=(X_{1},X_{2},X_{3})^{T}\sim N_{3}(\boldsymbol{\mu},\mathbf{\Sigma})$  with
\begin{align*}
\boldsymbol{\mu}=\left(\begin{array}{ccccccccccc}
1.2\\
0.7\\
3.0
\end{array}
\right)~and~
\mathbf{\Sigma}=\left(\begin{array}{ccccccccccc}
1.33&-0.067&2.63\\
-0.067&0.25&-0.50\\
2.63&-0.50&5.76
\end{array}
\right).
\end{align*}
Now, let $p=0.05$, $q=0.1,~0.2,~0.5,~0.8,~0.9,~0.95$, and $q=0.05$, \\$p=0.05,~0.1,~0.2,~0.5,~0.8,~0.9,$ then the $\mathrm{DTE}_{p,q}(X_{1})$ and $\mathrm{DTV}_{p,q}(X_{1})$ are presented in Figures 1 and 2, respectively:

From Figures 1 and 2, we can observe that when fixing $p=0.05$, DTE of $X_{1}$ is increasing with the increase of $q$, while DTV of $X_{1}$ is decreasing first and then is increasing with the increase of $q$; DTV of $X_{1}$ for $q=0.50$ is greatest; When fixing $q=0.95$, DTE and DTV of $X_{1}$ have the same trend in Figure 1.

Next, let $(p,q)=(0.05,0.75)$, $(p,q)=(0.10,0.80)$, $(p,q)=(0.15,0.85)$, $(p,q)=(0.20,0.90)$ and $(p,q)=(0.25,0.95)$, and let $(p,q)=(0.05,0.95)$, $(p,q)=(0.10,0.90)$, $(p,q)=(0.15,0.85)$, $(p,q)=(0.20,0.80)$ and $(p,q)=(0.25,0.75)$, then the $\mathrm{DTE}_{p,q}(X_{1})$ and $\mathrm{DTV}_{p,q}(X_{1})$ are showed in Figures 3 and 4, respectively:

From Figures 3 and 4, we can observe that when fixing $q-p=0.7$, DTE of $X_{1}$ is decreasing fist and then is decreasing with the increase of $p$, while DTV of $X_{1}$ is decreasing first, then is increasing and finally is decreasing with the increase of $p$; DTE of $X_{1}$ for $(p,q)=(0.10,0.80)$ is least; When fixing $p+q=1$, DTE of $X_{1}$ is the same, and is equal to $\mu_{1}=1.2$, while DTV of $X_{1}$ is decreasing with the increase of $p$.

Finally, let $\boldsymbol{p_{1}}=(0.10,~0.10,~0.10)^{T}$,~$\boldsymbol{q_{1}}=(0.80,~0.80,~0.80)^{T},$\\ $\boldsymbol{p_{2}}=(0.15,~0.15,~0.15)^{T}$,~$\boldsymbol{q_{2}}=(0.85,~0.85,~0.85)^{T},$ $\boldsymbol{p_{3}}=(0.20,~0.20,~0.20)^{T}$ and~$\boldsymbol{q_{3}}=(0.90,~0.90,~0.90)^{T},$ we obtain MDTEs of $\mathbf{X}$ presented in Table 1 and Figure 5:
\begin{table}[!htbp]
\centering
\begin{tabular}{|c||c|c|c|c|c|c|}
  \hline
  \backslashbox {$(\boldsymbol{p_{i}},\boldsymbol{q_{i}})$}{$MDTE_{(\boldsymbol{p_{i}},\boldsymbol{q_{i}})}(\mathbf{X})$}{$X_{i}$} &$X_{1}$& $X_{2}$&$X_{3}$\\
  \hline
  $(\boldsymbol{p_{1}},\boldsymbol{q_{1}})$ &1.027895 & 0.708670& 2.659672\\
  \hline
  $(\boldsymbol{p_{2}},\boldsymbol{q_{2}})$&1.200000 & 0.700000&  3.000000\\
  \hline
  $(\boldsymbol{p_{3}},\boldsymbol{q_{3}})$&1.372105 & 0.691330&  3.340328\\
  \hline
\end{tabular}\\
Table 1: The MDTE of $\mathbf{X}$ for $(\boldsymbol{p_{1}},\boldsymbol{q_{1}}),~(\boldsymbol{p_{2}},\boldsymbol{q_{2}})~ and ~(\boldsymbol{p_{3}},\boldsymbol{q_{3}})$.
\end{table}

From Table 1 and Figure 5, we find that for same $X_{k},~i=1,~2,~3$, the MDTE of $X_{1}$ and $X_{3}$ for $(\boldsymbol{p_{3}},\boldsymbol{q_{3}})$ are the greatest, while the MDTE of $X_{2}$ for $(\boldsymbol{p_{3}},\boldsymbol{q_{3}})$ is the least; the MDTE of $X_{1}$ and $X_{3}$ for $(\boldsymbol{p_{1}},\boldsymbol{q_{1}})$ are the least, while the MDTE of $X_{2}$ for $(\boldsymbol{p_{1}},\boldsymbol{q_{1}})$ is the greatest; For same $(\boldsymbol{p_{k}},\boldsymbol{q_{k}})$,~$k=1,~2,~3$, the MDTE of $X_{3}$ is the greatest, and the MDTE of $X_{2}$ is the least.

The MDTCovs of $\mathbf{X}$ are showed as follows:
\begin{align*}
\mathrm{MDTCov}_{(\boldsymbol{p_{1}},\boldsymbol{q_{1}})}(\mathbf{X})=\left(\begin{array}{ccccccccccc}
0.425369&-0.021428&0.841143\\
-0.021428&0.247704&-0.409885\\
0.841143&-0.409885&2.222636
\end{array}
\right),
\end{align*}
\begin{align*}
\mathrm{MDTCov}_{(\boldsymbol{p_{2}},\boldsymbol{q_{2}})}(\mathbf{X})=\left(\begin{array}{ccccccccccc}
0.411717&-0.020741&0.814146\\
-0.020741&0.247670&-0.408525\\
0.814146&-0.408525&2.169252
\end{array}
\right)
\end{align*}
and
\begin{align*}
\mathrm{MDTCov}_{(\boldsymbol{p_{3}},\boldsymbol{q_{3}})}(\mathbf{X})=\left(\begin{array}{ccccccccccc}
0.425369&-0.021428&0.841143\\
-0.021428&0.247704&-0.409885\\
0.841143&-0.409885&2.222636
\end{array}
\right).
\end{align*}
From $\mathrm{MDTCov}_{(\boldsymbol{p_{1}},\boldsymbol{q_{1}})}(\mathbf{X})$, $\mathrm{MDTCov}_{(\boldsymbol{p_{2}},\boldsymbol{q_{2}})}(\mathbf{X})$ and $\mathrm{MDTCov}_{(\boldsymbol{p_{3}},\boldsymbol{q_{3}})}(\mathbf{X})$, note that the main diagonal MDTCov for $(\boldsymbol{p_{2}},\boldsymbol{q_{2}})$ are smaller than that for corresponding others, and the MDTCov for $(\boldsymbol{p_{1}},\boldsymbol{q_{1}})$ and $(\boldsymbol{p_{3}},\boldsymbol{q_{3}})$ is the same.

 In addition, using $\mathrm{MDTCov}_{(\boldsymbol{p_{1}},\boldsymbol{q_{1}})}(\mathbf{X})$, $\mathrm{MDTCov}_{(\boldsymbol{p_{2}},\boldsymbol{q_{2}})}(\mathbf{X})$ and\\ $\mathrm{MDTCov}_{(\boldsymbol{p_{3}},\boldsymbol{q_{3}})}(\mathbf{X})$ and Eq. (\ref{(a9)}) we can give their MDTCorr matries, respectively:
\begin{align*}
\mathrm{MDTCorr}_{(\boldsymbol{p_{1}},\boldsymbol{q_{1}})}(\mathbf{X})=\left(\begin{array}{ccccccccccc}
1&-0.066013&0.865073\\
-0.066013&1&-0.552410\\
0.865073&-0.552410&1
\end{array}
\right),
\end{align*}
\begin{align*}
\mathrm{MDTCorr}_{(\boldsymbol{p_{2}},\boldsymbol{q_{2}})}(\mathbf{X})=\left(\begin{array}{ccccccccccc}
1&-0.064952&0.861485\\
-0.064952&1&-0.557349\\
0.861485&-0.557349&1
\end{array}
\right)
\end{align*}
and
\begin{align*}
\mathrm{MDTCorr}_{(\boldsymbol{p_{3}},\boldsymbol{q_{3}})}(\mathbf{X})=\left(\begin{array}{ccccccccccc}
1&-0.066013&0.865073\\
-0.066013&1&-0.552410\\
0.865073&-0.552410&1
\end{array}
\right).
\end{align*}

\section{Illustrative example}
We discuss MDTE and MDTCov of three industry segments in finance
(Banks $X_{1}$, Insurance $X_{2}$, Financial and Credit Service $X_{3}$) collecting stock return data in London stock exchange from April 2013 to November 2019 in the finance sector
of the market by using the results of parameter estimates in Shushi and Yao (2020). Using multivariate normal distribution to fit data. We denote it by $\mathbf{X}=(X_{1},X_{2},X_{3})^{T}\sim N_{3}(\boldsymbol{\mu},\mathbf{\Sigma})$. Parameters are computed using maximum likelihood estimation:
\begin{align*}
 \boldsymbol{\mu}=10^{-3}\left(\begin{array}{ccccccccccc}
-1.140677\\
5.896240\\
2.107343
\end{array}
\right),
\mathbf{\Sigma}=10^{-4}\left(\begin{array}{ccccccccccc}
19.088935&12.503116&-3.720492\\
12.503116&20.268816&-3.162601\\
-3.720492&-3.162601&8.851913
\end{array}
\right).
\end{align*}
Let $\boldsymbol{p_{1}}=(0.10,~0.10,~0.10)^{T}$,~$\boldsymbol{q_{1}}=(0.80,~0.80,~0.80)^{T},$ $\boldsymbol{p_{2}}=(0.15,~0.15,~0.15)^{T}$,\\~$\boldsymbol{q_{2}}=(0.85,~0.85,~0.85)^{T},$ $\boldsymbol{p_{3}}=(0.20,~0.20,~0.20)^{T}$ and~$\boldsymbol{q_{3}}=(0.90,~0.90,~0.90)^{T},$ then the MDTE of $\mathbf{X}$ is showed in Table 2 and Figure 6:
\begin{table}[!htbp]
\centering
\begin{tabular}{|c||c|c|c|c|c|c|}
  \hline
  \backslashbox {$(\boldsymbol{p_{i}},\boldsymbol{q_{i}})$}{$\mathrm{MDTE}_{(\boldsymbol{p_{i}},\boldsymbol{q_{i}})}(\mathbf{X})$}{$X_{i}$} &$X_{1}(10^{-3})$& $X_{2}(10^{-3})$&$X_{3}(10^{-3})$\\
  \hline
  $(\boldsymbol{p_{1}},\boldsymbol{q_{1}})$ &-7.660832 & 0.812427& 3.426994\\
  \hline
  $(\boldsymbol{p_{2}},\boldsymbol{q_{2}})$&-1.140677 & 5.896240&  2.107343\\
  \hline
  $(\boldsymbol{p_{3}},\boldsymbol{q_{3}})$&5.379478 & 10.980053&  0.787692\\
  \hline
\end{tabular}\\
Table 2: The MDTE of $\mathbf{X}$ for $(\boldsymbol{p_{1}},\boldsymbol{q_{1}}),~(\boldsymbol{p_{2}},\boldsymbol{q_{2}})~ and ~(\boldsymbol{p_{3}},\boldsymbol{q_{3}})$.
\end{table}

From Table 2 and Figure 6, we see that for same industry segments, the MDTE of Banks and Insurance for $(\boldsymbol{p_{3}},\boldsymbol{q_{3}})$ are the greatest, while the MDTE of Financial and Credit Service  for $(\boldsymbol{p_{3}},\boldsymbol{q_{3}})$ is the least; the MDTE of Banks and Insurance for $(\boldsymbol{p_{1}},\boldsymbol{q_{1}})$ are the least, while the MDTE of Financial and Credit Service for $(\boldsymbol{p_{1}},\boldsymbol{q_{1}})$ is the greatest. For different industry segments, there are different results for choosing different $(\boldsymbol{p},\boldsymbol{q})$. When $(\boldsymbol{p},\boldsymbol{q})=(\boldsymbol{p_{1}},\boldsymbol{q_{1}})$, the MDTE of Financial and Credit Service is the greatest, and the MDTE of Banks is the least, which means that there are the largest difference of MDTE between Banks and Financial and Credit Service; When $(\boldsymbol{p},\boldsymbol{q})=(\boldsymbol{p_{2}},\boldsymbol{q_{2}})$, the MDTE of Insurance is the greatest, and the MDTE of Banks is the least, which means that there are the largest difference of MDTE between Banks and Insurance; When $(\boldsymbol{p},\boldsymbol{q})=(\boldsymbol{p_{3}},\boldsymbol{q_{3}})$, the MDTE of Insurance is the greatest, and the MDTE of Financial and Credit Service is the least, which means that there are the largest difference of MDTE between Insurance and Financial and Credit Service.

Now, the MDTCovs of $\mathbf{X}$ are presented as follows:
\begin{align*}
\mathrm{MDTCov}_{(\boldsymbol{p_{1}},\boldsymbol{q_{1}})}(\mathbf{X})=10^{-4}\left(\begin{array}{ccccccccccc}
6.105151&4.063519&-1.193799\\
4.063519&14.063369&-1.476991\\
-1.193799&-1.476991&8.357612
\end{array}
\right),
\end{align*}
\begin{align*}
\mathrm{MDTCov}_{(\boldsymbol{p_{2}},\boldsymbol{q_{2}})}(\mathbf{X})=10^{-4}\left(\begin{array}{ccccccccccc}
5.909184&3.870473&-1.151718\\
3.870473&14.221969&-1.456490\\
-1.151718&-1.456490&8.349834
\end{array}
\right)
\end{align*}
and
\begin{align*}
\mathrm{MDTCov}_{(\boldsymbol{p_{3}},\boldsymbol{q_{3}})}(\mathbf{X})=10^{-4}\left(\begin{array}{ccccccccccc}
6.105151&4.063520&-1.193799\\
4.063520&14.063369&-1.476991\\
-1.193799&-1.476991&8.357612
\end{array}
\right).
\end{align*}
From $\mathrm{MDTCov}_{(\boldsymbol{p_{1}},\boldsymbol{q_{1}})}(\mathbf{X})$, $\mathrm{MDTCov}_{(\boldsymbol{p_{2}},\boldsymbol{q_{2}})}(\mathbf{X})$ and $\mathrm{MDTCov}_{(\boldsymbol{p_{3}},\boldsymbol{q_{3}})}(\mathbf{X})$, note that the main diagonal MDTCov for $(\boldsymbol{p_{2}},\boldsymbol{q_{2}})$ are least than that for corresponding others.
 \section{Concluding remarks}
 In this paper, we have defined MDTE and MDTCov risk measures, which are extending of MTCE and MTCov risk measures in Landsman et al. (2018). The expression of MDTE and MDTCov for elliptical distributions  have been derived. There are some special cases, including normal, student-$t$, logistic, Laplace and Pearson type VII distributions. As an illustrative example, the MDTE and MDTCov of three industry segments'
(Banks, Insurance, Financial and Credit Service) stock return in London stock exchange  are discussed. Note that, in general, let $VaR_{\boldsymbol{v}}(\boldsymbol{X})=\boldsymbol{x_{v}},~\boldsymbol{v=p,q}$ in (\ref{(22)}) and (\ref{(36)}), we can obtain general formulas of the MDTE and MDTCov for elliptical distributions. We know $\mathbf{Y}|W\sim N_{p}\left(\boldsymbol{\mu+\tau} W,~W\mathbf{\Sigma}\right)$ from Roozegar et al. (2020), so the MDTE and MDTCov of $\mathbf{Y}|W$ are special cases of the our results.
 Furthermore, in Landsman and Shushi (2021), the authors provided expressions of MTCE and MTCov for logistic-elliptical distributions. It will, therefore, be of interest to extend the results established here to the
logistic-elliptical distributions.
\section*{Acknowledgments}
\noindent  The research was supported by the National Natural Science Foundation of China (No. 12071251, 11571198)
\section*{Conflicts of Interest}
\noindent The authors declare that they have no conflicts of interest.
\section*{References}
\bibliographystyle{model1-num-names}







\end{document}